\documentclass[11pt]{article}

\usepackage{amsmath}
\usepackage{amsfonts}
\usepackage{amssymb}
\usepackage{amsthm}

\newtheorem{theorem}{Theorem}[section]
\newtheorem{proposition}[theorem]{Proposition}
\newtheorem{corollary}[theorem]{Corollary}
\newtheorem{lemma}[theorem]{Lemma}
\newtheorem{definition}[theorem]{Definition}
\newtheorem{notation}[theorem]{Notation}
\newtheorem{remark}[theorem]{Remark}

\newcommand{\cA}{ {\cal A} }

\newcommand{\bB}{ {\mathbb B} }
\newcommand{\bC}{ {\mathbb C} }
\newcommand{\cD}{ {\cal D} }

\newcommand{\cH}{ {\cal H}}

\newcommand{\cK}{ {\cal K}}
\newcommand{\cM}{ {\cal M} }

\newcommand{\bR}{ {\mathbb R} }

\newcommand{\cS}{ {\cal S} }

\newcommand{\bZ}{ {\mathbb Z} }

\newcommand{\leqleq}{ \ll }
\newcommand{\geqgeq}{ \gg }

\newcommand{\Int}{ \mbox{Int} } 
\newcommand{\cf}{ \mbox{Cf} } 
\newcommand{\Reta}{ \mbox{Reta} } 
\newcommand{\card}{ \mbox{card} } 
\newcommand{\term}{ \mbox{term} } 
\newcommand{\tleft}{ \triangleleft }

\newcommand{\Dalg}{ \cD_{\mathrm{alg}} }
\newcommand{\ncpolk}{ \bC \langle X_1, \ldots , X_k \rangle }
\newcommand{\ncserk}{ \bC_0 \langle \langle z_1, \ldots ,z_k \rangle \rangle }

\newcommand{\ecpi}{ \stackrel{\pi}{\sim} }

\newcommand{\freeplus}{ \boxplus }
\newcommand{\freetimes}{ \boxtimes }
\newcommand{\freestar}{ \framebox[7pt]{$\star$} }
\newcommand{\Uplus}{ \uplus }

\begin{document}

\title{\bf Free Brownian motion and evolution towards 
\boldmath{$\boxplus$}-infinite divisibility for \boldmath{$k$}-tuples }
\author{
Serban T. Belinschi 
\and Alexandru Nica \thanks{Research supported by a Discovery Grant of 
NSERC, Canada.} }

\date{ }

\maketitle

\begin{abstract}
Let $\cD_c (k)$ be the space of (non-commutative) distributions of 
$k$-tuples of selfadjoint elements in a $C^*$-probability space.
For every $t \geq 0$ we consider the transformation 
$\bB_t : \cD_c (k) \to \cD_c (k)$ defined by 
\[
{\mathbb B}_t ( \mu ) = \Bigl( \, \mu^{\boxplus (1+t)} \, 
\Bigr)^{\uplus (1/(1+t))}, \ \ \mu \in \cD_c (k),
\]
where $\boxplus$ and $\uplus$ are the operations of free additive 
convolution and respectively of Boolean convolution on $\cD_c (k)$.
We prove that 
$\bB_s \circ \bB_t = \bB_{s+t}$, $\forall \, s,t \geq 0$.
For $t=1$ we prove that $\bB_1 ( \, \cD_c (k) \, )$ is precisely 
the set $\cD_c^{\textrm{inf-div}} (k)$ of distributions in 
$\cD_c (k)$ which are infinitely divisible with respect to 
$\boxplus$, and that the map $\cD_c (k) \ni \mu 
\mapsto \bB_1 ( \mu ) \in \cD_c (k)^{\textrm{inf-div}}$
coincides with the multi-variable Boolean Bercovici-Pata bijection 
put into evidence in our previous paper \cite{BN06}. Thus for a 
fixed $\mu \in \cD_c (k)$, the process 
$\{ \bB_t ( \mu ) \mid t \geq 0 \}$ can be viewed as some kind of 
``evolution towards $\boxplus$-infinite divisibility''.

On the other hand we put into evidence a relation between the 
transformations $\bB_t$ and free Brownian motion. More precisely,
we introduce a map $\Phi : \cD_c (k) \to \cD_c (k)$ which 
transforms the free Brownian motion started at an arbitrary 
$\nu \in \cD_c (k)$ into the process 
$\{ \bB_t ( \mu ) \mid t \geq 0 \}$ for $\mu = \Phi ( \nu )$.
\end{abstract}

$\ $

\setcounter{section}{1}
\begin{center}
{\large\bf 1. Introduction}
\end{center}

{\bf 1.1 Review of past work.}
The study of noncommutative forms of independence for random 
variables has led to several ``convolution operations'' that 
can be defined on the space $\cM$ of probability distributions 
on $\bR$. Two such operations are the free (additive) convolution
$\boxplus$ and the Boolean convolution $\uplus$, which reflect 
the operations of addition of freely independent and 
respectively Boolean independent random variables.

In the paper \cite{BN07} we introduced a family 
$\{ \bB_t \mid t \geq 0 \}$ of transformations of $\cM$, defined 
by the formula
\begin{equation}  \label{eqn:1.1}
{\mathbb B}_t ( \mu ) = \Bigl( \, \mu^{\boxplus (1+t)} \, 
\Bigr)^{\uplus (1/(1+t))}, \ \ \forall \, t \geq 0,
\ \forall \, \mu \in \cM.
\end{equation}
The transformations $\bB_{t}$ turn out to form a semigroup:
$\bB_s \circ \bB_t = \bB_{s+t}, \ \  \forall \, s,t \geq 0$.
On the other hand for $t=1$ it turns out that the range set 
$\bB_1 ( \cM )$ is precisely the set $\cM^{\textrm{inf-div}}$
of distributions in $\cM$ which are infinitely divisible with 
respect to $\boxplus$; and moreover, the map
$\cM \ni \mu \mapsto \bB_1 ( \mu ) \in \cM^{\textrm{inf-div}}$
coincides with a remarkable bijection discovered by Bercovici and 
Pata \cite{BP99} in their study of relations between infinite 
divisibility with respect to $\boxplus$ and to $\uplus$. 

Due to the above properties of the transformations $\bB_t$,
for a fixed $\mu \in \cM$ the process $t \mapsto \bB_t ( \mu )$
can be viewed as a kind of ``evolution towards $\boxplus$-infinite
divisibility'' (where infinite divisibility is always reached by 
the time $t=1$). In \cite{BN07} it was observed that this process 
is related to free Brownian motion. Recall that the free 
Brownian motion started at $\nu \in \cM$ is the process 
$\{ \nu \boxplus \gamma_t \mid t \geq 0 \}$, where 
$\gamma_t \in \cM$ is the centered semicircular distribution of 
variance $t$. The connection between this and the transformations
$\bB_t$ is described as follows.

For a distribution $\mu \in \cM$ let $G_{\mu}$ and $F_{\mu}$
denote the Cauchy transform and respectively the reciprocal 
Cauchy transform of $\mu$; that is, we have
\[
G_{\mu} (z)= \int_{\bR} \frac{d \mu (s)}{z-s},
\mbox{ and } 
F_{\mu} (z)= 1/ G_{\mu} (z), \ \ z \in \bC \setminus \bR.
\]
By using basic facts from the theory of the Cauchy transform,
one easily sees that for every distribution $\nu \in \cM$ there 
exists a unique $\mu \in \cM$ such that 
\begin{equation}  \label{eqn:1.2}
F_{\mu} (z) = z - G_{\nu} (z), \ \ 
\forall \, z \in \bC \setminus \bR .
\end{equation}
One can thus define a map $\Phi : \cM \to \cM$ by putting 
$\Phi ( \nu ) := \mu$ with $\mu$ and $\nu$ related as in 
(\ref{eqn:1.2}). The map $\Phi$ turns out to be one-to-one, with 
image $\Phi ( \cM )$ consisting precisely of those distributions
$\mu \in \cM$ which have
$\int_{- \infty}^{\infty} t^2 \ d \mu (t) = 1$ and
$\int_{- \infty}^{\infty} t \ d \mu (t) = 0$. (A detailed 
presentation of these facts appears in Section 2 of the paper
\cite{M92}.) The relation between the transformations $\bB_t$ 
and free Brownian motion can be expressed by using the map 
$\Phi$ and the following formula:
\begin{equation}  \label{eqn:1.3}
\Phi ( \, \nu \boxplus \gamma_t \, ) = 
\bB_t ( \, \Phi ( \nu ) \, ), \ \ \forall \, \nu \in \cM,
\ \forall \, t > 0.
\end{equation}   
In other words the free Brownian motion started at $\nu$ 
corresponds exactly, via $\Phi$, to the process 
$\{ \bB_t ( \mu ) \mid t \geq 0 \}$ started at $\mu = \Phi ( \nu )$.

$\ $

{\bf 1.2 Description of results of this paper.}
In this paper we find multi-variable analogues for the results 
described above. Let $k$ be a positive integer, and let $\cD_c (k)$
denote the space of non-commutative distributions of $k$-tuples 
of selfadjoint elements in a $C^*$-probability space. The 
convolution operations $\boxplus$ and $\uplus$ can be defined on
$\cD_c (k)$, and for every $\mu \in \cD_c (k)$ it makes sense to
define convolution powers $\mu^{\boxplus p}$, $\forall \, p \geq 1$ 
and $\mu^{\uplus q}$, $\forall \, q > 0$. One can thus define a 
family $\{ \bB_t \mid t \geq 0 \}$ of transformations of 
$\cD_c (k)$ by exactly the same formula as in (\ref{eqn:1.1}):
\begin{equation}  \label{eqn:1.4}
{\mathbb B}_t ( \mu ) = \Bigl( \, \mu^{\boxplus (1+t)} \, 
\Bigr)^{\uplus (1/(1+t))}, \ \ \forall \, t \geq 0,
\ \forall \, \mu \in \cD_c (k).
\end{equation}
We prove that $\bB_s \circ \bB_t = \bB_{s+t}$, 
$\forall \, s,t \geq 0$. For $t=1$ we prove that 
$\bB_1 ( \, \cD_c (k) \, )$ is precisely the set 
$\cD_c^{\textrm{inf-div}} (k)$ of distributions in $\cD_c (k)$ 
which are infinitely divisible with respect to $\boxplus$, and
that the map $\cD_c (k) \ni \mu 
\mapsto \bB_1 ( \mu ) \in \cD_c (k)^{\textrm{inf-div}}$
coincides with the multi-variable Boolean Bercovici-Pata bijection 
put into evidence in our previous paper \cite{BN06}. Thus for a 
fixed $\mu \in \cD_c (k)$, the process 
$\{ \bB_t ( \mu ) \mid t \geq 0 \}$
can still be viewed as a kind of evolution towards 
$\boxplus$-infinite divisibility, which is now taking place in 
the framework of $\cD_c (k)$.

Moreover, we prove that the transformations $\bB_t$ relate to 
the multi-variable free Brownian motion in a similar way to the 
one presented above in the 1-dimensional case. The free Brownian 
motion started at $\nu \in \cD_c (k)$ is the process 
$\{ \nu \boxplus \gamma_t \mid t \geq 0 \}$ where 
$\gamma_t \in \cD_c (k)$ now stands for the joint distribution 
of a free family $x_1, \ldots , x_k$ of selfadjoint elements in 
a $C^*$-probability space, such that every $x_i$ has a centered 
semicircular distribution of variance $t$. In order to connect this
to the transformations $\bB_t$, we use a multi-variable analogue
for the map $\Phi$ which had been defined via Equation 
(\ref{eqn:1.2}). The multi-variable version of Equation 
(\ref{eqn:1.2}) involves formal power series in $k$ non-commuting
indeterminates $z_1, \ldots , z_k$ (instead of complex analytic
functions of one variable $z$). For $\mu \in \cD_c (k)$ let 
$M_{\mu}$ be its moment series,
\[
M_{\mu} (z_1, \ldots , z_k) := \sum_{n=1}^{\infty} \
\sum_{i_1, \ldots , i_n =1}^k \ \mu (X_{i_1} \cdots X_{i_n})
\, z_{i_1} \cdots z_{i_n};
\]
and let us moreover denote
\[
\eta_{\mu} := M_{\mu} \, (1+ M_{\mu})^{-1}, \ \ 
\mu \in \cD_c (k),
\]
where $(1+ M_{\mu})^{-1}$ is the inverse of $1+M_{\mu}$ under
multiplication in the ring of power series 
$\bC \langle \langle z_1, \ldots , z_k \rangle \rangle$. 
With these notations, the multi-variable 
analogue for (\ref{eqn:1.2}) is 
\begin{equation}  \label{eqn:1.5}
\eta_{\mu} (z_1, \ldots , z_k) = \sum_{i=1}^k \, z_i
\Bigl( 1+ M_{\nu} (z_1, \ldots , z_k) \Bigr) z_i
\end{equation}
(equality of formal power series). We prove that for every 
$\nu \in \cD_c (k)$ there exists a unique $\mu \in \cD_c (k)$ 
such that (\ref{eqn:1.5}) holds. One can thus define a map 
$\Phi : \cD_c (k) \to \cD_c (k)$ by putting $\Phi ( \nu ) := \mu$
with $\mu, \nu$ as in (\ref{eqn:1.5}), and it turns out that 
we then have the analogue of (\ref{eqn:1.3}):
\begin{equation}  \label{eqn:1.6}
\Phi ( \, \nu \boxplus \gamma_t \, ) = 
\bB_t ( \, \Phi ( \nu ) \, ), \ \ \forall \, \nu \in \cD_c (k),
\ \forall \, t > 0.
\end{equation}   

Concerning the seemingly different appearance of Equations 
(\ref{eqn:1.2}) and (\ref{eqn:1.5}), we make the following comment.
Suppose that $k=1$ and that $\mu , \nu \in \cD_c (1)$ are 
identified as compactly supported distributions on $\bR$. Then
$\eta_{\mu}$ and $M_{\nu}$ can be viewed as analytic functions on
a neighbourhood of 0, and (for $z$ running in a suitable domain of 
$\bC$) we have
\[
F_{\mu} (z) = z \Bigl( \, 1 - \eta_{\mu} ( 1/z ) \, \Bigr) , \ \ 
G_{\nu} (z) = \frac{1}{z} \Bigl( \, 1+ M_{\nu} 
( 1/z ) \, \Bigr) .
\]
By substituting these formulas into Equation (\ref{eqn:1.2}), 
and by replacing $z$ with $1/z$, we bring (\ref{eqn:1.2}) to the 
form
\[
\eta_{\mu} (z) = z^2 \Bigl( 1 + M_{\nu} (z) \Bigr),
\]
which is exactly the 1-dimensional version of Equation (\ref{eqn:1.5}).

$\ $

{\bf 1.3 Further remarks.}
The results in \cite{BN07} were proved by using complex analytic 
functions. The methods used in this paper are completely different, 
they rely on the combinatorics of non-crossing partitions. 
For most of the paper we use the larger algebraic framework of 
$\Dalg (k)$, the space of all possible joint distributions of 
$k$-tuples in a (sheer algebraic) non-commutative probability 
space. It makes sense to define $\bB_t$ as a bijective transformation 
of $\Dalg (k)$, then prove the algebraic statements made about 
$\{ \bB_t \mid t \geq 0 \}$ in this larger framework; these
properties of the $\bB_t$ are summarized in Theorem 
\ref{thm:4.11} of the paper. In the same theorem we also point out 
an additional property of $\bB_t$, that
\begin{equation}  \label{eqn:1.7}
\bB_t ( \mu \boxtimes \nu ) = \bB_t ( \mu ) \boxtimes \bB_t ( \nu ),
\ \ \forall \, t \geq 0, \ \forall \, \mu , \nu \in \Dalg (k),
\end{equation}
where $\boxtimes$ (the free {\em multiplicative} convolution) 
is the operation on $\Dalg (k)$ which corresponds to the 
multiplication of free $k$-tuples of random variables in a 
non-commutative probability space.

The map $\Phi$ defined via Equation (\ref{eqn:1.5}) and the 
relation between the transformations $\bB_t$ and free Brownian 
motion given in (\ref{eqn:1.6}) can be considered on
$\Dalg (k)$ as well. The proof of formula (\ref{eqn:1.6}) is in
fact done in this algebraic framework, in Theorem \ref{thm:6.2} 
of the paper.

After the algebraic results are established, what remains to be
done is make sure that the transformations $\bB_t$ do indeed the 
job they are supposed to, when they are restricted to the smaller
space $\cD_c (k)$. Some of the verifications needed here are direct
consequences of things proved in our preceding paper \cite{BN06}.
But there is one verification that we are left with,
namely that the map $\Phi : \Dalg (k) \to \Dalg (k)$ carries 
$\cD_c (k)$ into itself. We prove this fact by providing an
{\em operator model} for how $\Phi$ works on $\cD_c (k)$; this 
operator model is discussed in Remark \ref{rem:7.3} and 
Theorem \ref{thm:7.4} of the paper.

We conclude the introduction with an outline of how the paper is 
organized. Throughout the whole paper $k$ is a fixed positive
integer -- the ``number of indeterminates'' we are working 
with. In Section 2 we review the algebraic framework of
$\Dalg (k)$, and the $R$ and $\eta$ transforms for distributions 
in $\Dalg (k)$. Section 3 is a review section as well, devoted to 
an important bijection on power series introduced in \cite{BN06}, 
the bijection ``$\Reta$'' sending $R_{\mu} \mapsto \eta_{\mu}$ 
for every $\mu \in \Dalg (k)$); this bijection is the workhorse 
for many of the computations with power series done in the present 
paper. In Section 4 we introduce $\bB_t$ as a bijective 
transformation of $\Dalg (k)$ and we prove some of its basic 
properties; the results of the section are summarized in Theorem 
\ref{thm:4.11}. In Section 5 we establish a formula for moments 
of the free Brownian motion, which is needed in the proof of the 
connection between free Brownian motion and the $\bB_t$. The proof 
of this connection is then done in Section 6, Theorem \ref{thm:6.2}. 
The final Section 7 of the paper deals with the framework of 
$\cD_c (k)$, the main point 
of the section being the operator model for $\Phi$.

$\ $

$\ $

\setcounter{section}{2}
\begin{center}
{\large\bf 2. Non-commutative convolutions and transforms on 
              \boldmath{$\Dalg (k)$} }
\end{center}

\setcounter{equation}{0}
\setcounter{theorem}{0}

\begin{definition}  \label{def:2.1}
(Non-commutative distributions.)

{\rm
$1^o$ We denote by $\ncpolk$ the algebra of non-commutative 
polynomials in $X_1, \ldots ,$
$X_k$. Thus $\ncpolk$ has a linear basis
\begin{equation}  \label{eqn:3.101}
\{ 1 \} \cup \{ X_{i_1} \cdots X_{i_n} \mid n \geq 1, \
1 \leq i_1, \ldots , i_n \leq k \},
\end{equation}
where the monomials in the basis are multiplied by concatenation.
When needed, $\bC \langle X_1, \ldots ,$
$X_k \rangle$ will be viewed as a $*$-algebra, with 
$*$-operation determined uniquely by the fact that each of 
$X_1, \ldots , X_k$ is selfadjoint.

$2^o$ Let $( \cA , \varphi )$ be a non-commutative probability 
space; that is, $\cA$ is a unital algebra over $\bC$, and 
$\varphi : \cA \to \bC$ is a linear functional, normalized by 
the condition that $\varphi ( 1_{\cA} ) = 1$. 
For $x_1, \ldots , x_k \in \cA$, the {\em joint distribution}
of $x_1, \ldots , x_k$ is the linear functional 
$\mu_{x_1, \ldots , x_k} : \ncpolk \to \bC$ which acts on the 
linear basis (\ref{eqn:3.101}) by the formula
\begin{equation}  \label{eqn:3.102}
\left\{  \begin{array}{lcl}
\mu_{x_1, \ldots , x_k} ( 1 ) & = & 1                        \\
\mu_{x_1, \ldots , x_k} ( X_{i_1} \cdots X_{i_n} ) & = & 
                           \varphi ( x_{i_1} \cdots x_{i_n} ),     \\
      & & \forall \, n \geq 1, \ 1 \leq i_1, \ldots , i_n \leq k. 
\end{array}   \right.
\end{equation}

$3^o$ We will denote
\begin{equation}  \label{eqn:4.303}
\Dalg (k) := \{ \mu : \ncpolk \to \bC \mid \mu \mbox{ linear, }
\mu ( 1 ) = 1 \} .
\end{equation}
It is immediate that $\Dalg (k)$ is precisely the set of linear 
functionals on $\ncpolk$ that can arise as joint distribution for 
some $k$-tuple $x_1, \ldots , x_k$ in a non-commutative probability 
space. }
\end{definition}

\vspace{6pt}

\begin{remark} \label{rem:2.2}
(The operations $\boxplus$ and $\uplus$ on $\Dalg (k)$.)
{\rm
These operations are defined via the general principle that if one
has a form of independence for non-commutative random variables, 
then the addition of independent $k$-tuples of random variables 
will induce a ``convolution'' operation on $\Dalg (k)$. 

The operation $\boxplus$ arises in this way, in connection to the 
concept of {\em free independence}. Given $\mu , \nu \in \Dalg (k)$,
one can always find random variables 
$x_1, \ldots , x_k, y_1, \ldots , y_k$ in a non-commutative 
probability space $( \cA, \varphi )$ such that the joint 
distribution of the $k$-tuple 
$x_1, \ldots , x_k$ is equal to $\mu$, the joint distribution of 
$y_1, \ldots , y_k$ is equal to $\nu$, and such that 
$\{ x_1, \ldots , x_k \}$ is freely independent from 
$\{ y_1, \ldots , y_k \}$ in $( \cA , \varphi )$. The joint 
distribution of the $k$-tuple $x_1 + y_1, \ldots , x_k +y_k$ turns 
out to depend only on $\mu$ and $\nu$; and the free additive 
convolution $\mu \boxplus \nu$ is equal, by definition, to the joint 
distribution of $x_1 + y_1, \ldots , x_k + y_k$. 

The operation $\uplus$ is defined in the same way, but where we 
use the concept of {\em Boolean independence}: given 
$\mu , \nu \in \Dalg (k)$, the Boolean convolution $\mu \uplus \nu$ 
is the (uniquely determined) distribution of 
$x_1 + y_1, \ldots , x_k + y_k$, where 
the joint distribution of $x_1, \ldots , x_k$ is equal to $\mu$, 
the joint distribution of $y_1, \ldots , y_k$ is equal to $\nu$, and 
$\{ x_1, \ldots , x_k \}$ is Boolean independent from 
$\{ y_1, \ldots , y_k \}$.   

A commonly used method for studying the operations $\boxplus$ and 
$\uplus$ goes by considering {\em cumulants} for distributions in 
$\Dalg (k)$: in relation to $\boxplus$ one considers the 
{\em free cumulants} introduced in \cite{S94}, while for $\uplus$ 
one uses the {\em Boolean cumulants} which go back all the way 
to \cite{W73}. In this paper we will work with the concepts, 
equivalent to cumulants, of {\em linearizing transforms} for 
$\boxplus$ and $\uplus$. Specifically, for a distribution 
$\mu \in \Dalg (k)$ we will work with the {\em R-transform}
$R_{\mu}$ (a formal power series which records the free cumulants
of $\mu$) and with the the {\em $\eta$-series} $\eta_{\mu}$
(which does the same job in connection to the Boolean cumulants
of $\mu$). The precise definitions of $R_{\mu}$ and 
$\eta_{\mu}$ will be reviewed in the next notation and remark.
The meaning of the statement that $R$ and $\eta$ are 
``linearizing transforms'' for $\boxplus$ and respectively for 
$\uplus$ is that we have:
\begin{equation}   \label{eqn:2.021}
R_{\mu \freeplus \nu} = R_{\mu} + R_{\nu},
\ \ \forall \, \mu , \nu \in \Dalg (k),
\end{equation}
and
\begin{equation}   \label{eqn:2.022}
\eta_{\mu \Uplus \nu} = \eta_{\mu} + \eta_{\nu},
\ \ \forall \, \mu , \nu \in \Dalg (k).
\end{equation}  }
\end{remark}

$\ $

\begin{notation} \label{def:2.3}
{\rm
We will denote by $\ncserk$
the space of power series with complex coefficients and with 
vanishing constant term, in $k$ non-commuting indeterminates 
$z_1, \ldots , z_k$. The general form of a series 
$f \in \ncserk$ is thus
\begin{equation}  \label{eqn:5.101}
f( z_1, \ldots , z_k) = \sum_{n=1}^{\infty} \
\sum_{i_1, \ldots , i_n =1}^k  \ \alpha_{(i_1, \ldots , i_n)}
z_{i_1} \cdots z_{i_n},
\end{equation}
where the coefficients  $\alpha_{(i_1, \ldots , i_n)}$ are from $\bC$. }
\end{notation}

\vspace{6pt}

\begin{remark} \label{rem:2.4}
(review of the series $M_{\mu}$, $R_{\mu}$, $\eta_{\mu}$).
{\rm 
Let $\mu$ be a distribution in $\Dalg (k)$. 

$1^o$ We will denote by $M_{\mu}$ the series in $\ncserk$ given by
\begin{equation}  \label{eqn:2.041}
M_{\mu} (z_1, \ldots , z_k) := \sum_{n=1}^{\infty} \
\sum_{i_1, \ldots , i_n =1}^k \ \mu (X_{i_1} \cdots X_{i_n})
\, z_{i_1} \cdots z_{i_n}.
\end{equation}
$M_{\mu}$ is called the {\em moment series} of $\mu$, and
its coefficients ($\mu ( X_{i_1} \cdots X_{i_n} )$, with $n \geq 1$
and $1 \leq i_1, \ldots , i_n \leq k$) are called the 
{\em moments} of $\mu$.

$2^o$ The {\em $\eta$-series} of $\mu$ is
\begin{equation}  \label{eqn:2.042}
\eta_{\mu} := M_{\mu} (1+M_{\mu})^{-1} \in \ncserk ,
\end{equation}
where $(1+ M_{\mu})^{-1}$ is the inverse of $1 + M_{\mu}$ under
multiplication, in the ring of series 
$\bC \langle \langle z_1, \ldots , z_k \rangle \rangle$.
The coefficients of $\eta_{\mu}$ are called the
{\em Boolean cumulants} of $\mu$.

$3^o$ There exists a unique series $R_{\mu} \in \ncserk$ which 
satisfies the functional equation
\begin{equation}  \label{eqn:2.043}
R_{\mu} \Bigl( \, z_1 (1+M_{\mu}), \ldots , z_k (1+M_{\mu} \,
\Bigr) = M_{\mu}.
\end{equation}
Indeed, it is easily seen that Equation (\ref{eqn:2.043}) 
amounts to a recursion which determines uniquely the coefficients
of $R_{\mu}$ in terms of those of $M_{\mu}$. The series $R_{\mu}$
is called the {\em R-transform} of $\mu$, and its coefficients
are called the {\em free cumulants} of $\mu$. (See the discussion
in Lecture 16 of \cite{NS06}, and specifically Theorem 16.15 and
Corollary 16.16 of that lecture.) }
\end{remark}

\vspace{6pt}

\begin{remark}  \label{rem:2.5}
{\rm 
It is immediate that for every $f \in \ncserk$ there exists a 
unique distribution $\mu \in \Dalg (k)$ such that $\eta_{\mu} =f$.
(This is because, as one immediately checks, the equation 
$M_{\mu} ( 1+M_{\mu} )^{-1} = f$ is equivalent to 
$M_{\mu} = f(1-f)^{-1}$.) Thus the map $\mu \mapsto \eta_{\mu}$
is a bijection from $\Dalg (k)$ onto $\ncserk$.

Likewise, the map $\Dalg (k) \ni \mu \mapsto R_{\mu} \in \ncserk$
is bijective. The fact that for every $g \in \ncserk$ there exists
a unique $\mu \in \Dalg (k)$ such that $R_{\mu} = g$ is easily 
seen when one writes explicitly the relations between the 
coefficients of $R_{\mu}$ and $M_{\mu}$ that are coming out of
(\ref{eqn:2.043}) -- see Lectures 11 and 16 of \cite{NS06}.  }
\end{remark}

\vspace{6pt}

\begin{remark}  \label{rem:2.6}
(Convolution powers.)
{\rm For $\mu \in \Dalg (k)$ and a positive integer $n$ one 
denotes the $n$-fold convolution 
$\mu \boxplus \cdots  \boxplus \mu$ by $\mu^{\boxplus n}$. From 
the additivity (\ref{eqn:2.021}) of the $R$-transform 
it follows that $R_{\mu^{\boxplus n}} = n \cdot R_{\mu}$, 
and the latter formula can be extended to the case when $n$
is not necessarily integer. More precisely, for every 
$\mu \in \Dalg (k)$ and $t \in (0, \infty )$ one defines 
the convolution power $\mu^{\boxplus t}$ to be the unique 
distribution in $\Dalg (k)$ which has $R$-transform equal to
\begin{equation}  \label{eqn:2.061}
R_{\mu^{\boxplus t}} = t \cdot R_{\mu}.
\end{equation}
It is immediate that the $\boxplus$-convolution powers defined
in this way satisfy the usual rules for operating with exponents:
\begin{equation}  \label{eqn:2.062}
\mu^{\boxplus s} \boxplus \mu^{\boxplus t} =
\mu^{\boxplus (s+t)} \mbox{ and } 
\Bigl( \, \mu^{\boxplus s} \, \Bigr)^{\boxplus t} =
\mu^{\boxplus st}, \ \ \forall \, \mu \in \Dalg (k), \ 
\forall \, s,t > 0. 
\end{equation}
Note that, as a consequence, one has that for every fixed 
$t \in (0, \infty )$ the map $\mu \mapsto \mu^{\boxplus t}$ is 
a bijection from $\Dalg (k)$ onto itself.

A similar discussion can be made in connection to the convolution
powers with respect to $\uplus$: for every $\mu \in \Dalg (k)$ and 
$t \in ( 0, \infty )$ one defines the convolution power 
$\mu^{\uplus t}$ to be the unique distribution in $\Dalg (k)$ 
which has $\eta$-series equal to
\begin{equation}  \label{eqn:2.063}
\eta_{\mu^{\uplus t}} = t \cdot \eta_{\mu}.
\end{equation}
Then the $\uplus$-convolution powers satisfy the usual rules 
of operating with exponents, and for every fixed 
$t \in (0, \infty )$ the map  $\mu \mapsto \mu^{\uplus t}$ is 
a bijection from $\Dalg (k)$ onto itself.   } 
\end{remark}

$\ $

$\ $

\setcounter{section}{3}
\begin{center}
{\large\bf 3. The bijection ``Reta'', and its combinatorial properties}
\end{center}

\setcounter{equation}{0}
\setcounter{theorem}{0}

\begin{definition}  \label{def:3.1}
{\rm 
We will denote
\begin{equation}  \label{eqn:3.011}
\Reta := \eta \circ R^{-1} : \ncserk \to \ncserk ,
\end{equation}
where $R, \eta : \Dalg (k) \to \ncserk$ are the bijections 
$\mu \mapsto R_{\mu}$ and respectively $\mu \mapsto \eta_{\mu}$
that were discussed in Remark \ref{rem:2.5}. In other words,
$\Reta$ is the bijection from $\ncserk$ onto itself which is 
uniquely determined by the requirement that
\begin{equation}  \label{eqn:3.012}
\Reta ( R_{\mu}) = \eta_{\mu}, \ \ \forall \, \mu \in \Dalg (k).
\end{equation} }
\end{definition}

\vspace{6pt}

\begin{remark}  \label{rem:3.15}
{\rm The bijection $\Reta$ was introduced in our previous paper
\cite{BN06}. Its name was chosen by looking at Equation 
(\ref{eqn:3.012}) (the transformation of $\ncserk$
that ``converts $R$ into $\eta$''). It is very useful that one 
can alternatively describe $\Reta$ via an explicit formula which
gives directly the coefficients of the series $\Reta (f)$ in terms
of those of $f$, for $f \in \ncserk$. This formula is reviewed 
(following \cite{BN06}) in Proposition \ref{prop:3.4} below. 
It involves summations indexed by non-crossing partitions, and 
in order to present it we will start with a very concise review 
(intended mostly for setting notations) of the lattice $NC(n)$ of
non-crossing partitions. For a more detailed introduction to 
$NC(n)$ and to how it is used in free probability, we refer to 
\cite{NS06}, Lectures 9 and 10.  }
\end{remark}

\vspace{6pt}

\begin{remark}  \label{rem:3.2}
(Review of $NC(n)$.)
{\rm   Let $n$ be a positive integer.

$1^o$ Let $\pi$ = $\{ B_{1} , \ldots , B_{p} \}$ be a partition 
of $\{ 1, \ldots ,n \}$ -- i.e. $B_{1} , \ldots , B_{p}$ are 
pairwise disjoint non-void sets (called the {\em blocks} of $\pi$), 
and $B_{1} \cup \cdots \cup B_{p}$ = 
$\{ 1, \ldots , n \}$. We say that $\pi$ is {\em non-crossing} if for 
every $1 \leq i < j < k < l \leq n$ such that $i$ is in the same
block with $k$ and $j$ is in the same block with $l$, it necessarily
follows that all of $i,j,k,l$ are in the same block of $\pi$.
The set of all non-crossing partitions of 
$\{ 1, \ldots , n \}$ will be denoted by $NC(n).$ 

$2^o$ For $\pi \in NC(n)$, the number of blocks of $\pi$ will be 
denoted by $| \pi |$. 

$3^o$ On $NC(n)$ we consider the partial order given 
by {\em reversed refinement}: for $\pi , \rho \in NC(n)$, we write
``$\pi \leq \rho$'' to mean that every block of $\rho$ is a union of
blocks of $\pi$. The minimal and maximal element of $( NC(n), \leq )$ 
are denoted by $0_n$ (the partition of $\{ 1, \ldots , n \}$ into 
$n$ singleton blocks) and respectively $1_n$ (the partition of 
$\{ 1, \ldots , n \}$ into only one block).  

$4^o$ In the considerations about $\Reta$, an important role is 
played by another partial order relation on $NC(n)$, which was 
introduced in \cite{BN06} and is denoted by ``$\leqleq$''.
For $\pi , \rho \in NC(n)$ we will write
``$\pi \leqleq \rho$'' to mean that $\pi \leq \rho$ and that, in 
addition, the following condition is fulfilled:
\begin{equation}  
\left\{  \begin{array}{l}
\mbox{For every block $C$ of $\rho$ there exists a block}  \\
\mbox{$B$ of $\pi$ such that $\min (C), \max (C) \in B$.}
\end{array}  \right.
\end{equation}

It is immediately verified that ``$\leqleq$'' is indeed a partial 
order relation on $NC(n)$. It is much coarser than the reversed 
refinement order. For instance, the inequality $\pi \leqleq 1_n$ is not 
holding for all $\pi \in NC(n)$, but it rather amounts to the condition 
that the numbers $1$ and $n$ belong to the same block of $\pi$. 
At the other end of $NC(n)$, the inequality $\pi \geqgeq 0_n$ can 
only take place when $\pi = 0_n$.  }
\end{remark}

\vspace{6pt}

\begin{definition}  \label{def:3.3}
(coefficients for series in $\ncserk$).

{\rm 
$1^o$ For $n \geq 1$ and $1 \leq i_1, \ldots , i_n \leq k$ we will 
denote by
\begin{equation}  \label{eqn:3.031}
\cf_{(i_1, \ldots , i_n)} :
\bC_0 \langle \langle z_1, \ldots , z_k \rangle \rangle \to \bC
\end{equation}
the linear functional which extracts the coefficient of 
$z_{i_1} \cdots z_{i_n}$ in a series 
$f \in \bC_0 \langle \langle z_1, \ldots , z_k \rangle \rangle$.
Thus for $f$ written as in Equation (\ref{eqn:5.101}) we have 
$\cf_{(i_1, \ldots , i_n)} (f) = \alpha_{(i_1, \ldots , i_n)}$.

$2^o$ Suppose we are given a positive integer $n$, some indices 
$i_1, \ldots , i_n \in \{ 1, \ldots , k \}$, and a partition 
$\pi \in NC(n)$. We define a (generally non-linear) functional 
\begin{equation}  \label{eqn:3.032}
\cf_{(i_1, \ldots , i_n) ; \pi} :
\bC_0 \langle \langle z_1, \ldots , z_k \rangle \rangle \to \bC ,
\end{equation}
as follows. For every block $B = \{ b_1, \ldots , b_m \}$ of $\pi$, 
with $1 \leq b_1 < \cdots < b_m \leq n$, let us use the notation 
\[
(i_1, \ldots , i_n) \vert B \ := \ (i_{b_1}, \ldots , i_{b_m})
\in \{ 1, \ldots , k \}^m.
\]
Then we define
\begin{equation}  \label{eqn:5.102}
\cf_{(i_1, \ldots , i_n); \pi} (f) \ := \
\prod_{B \ \mathrm{block \ of} \ \pi} \ \cf_{(i_1, \ldots , i_n)|B} (f),
\ \ \forall \, f \in \ncserk .
\end{equation}
(For example if we had $n=5$ and $\pi = \{ \{ 1,4,5 \} , \{ 2,3 \} \}$, 
and if $i_1, \ldots , i_5$ would be some fixed indices from 
$\{ 1, \ldots , k \}$, then the above formula would become
\[
\cf_{(i_1, i_2, i_3, i_4, i_5) ; \pi } (f) \ = \
\cf_{(i_1, i_4, i_5)} (f) \cdot
\cf_{(i_2, i_3)} (f),
\]
$f \in \bC_0 \langle \langle z_1, \ldots , z_k \rangle \rangle$.) 
The quantities $\cf_{(i_1, \ldots , i_n); \pi} (f)$ will be 
referred to as {\em generalized coefficients} of the series $f$. }
\end{definition}

\vspace{6pt}

\begin{proposition}  \label{prop:3.4}
Let $f,g$ be series in $\ncserk$ such that
$\Reta (f) =g$. Then for every $n \geq 1$ and 
$1 \leq i_1, \ldots , i_n \leq k$ we have
\begin{equation}  \label{eqn:3.041}
\cf_{(i_1, \ldots , i_n)} (g) 
 = \sum_{ \begin{array}{c}
{\scriptstyle  \pi \in NC(n),}  \\
{\scriptstyle \pi \leqleq 1_n}
\end{array}  } \ \cf_{(i_1, \ldots , i_n); \pi} (f).
\end{equation}
More generally, we have the following formula for a
generalized coefficient $\cf_{(i_1, \ldots , i_n); \rho} (g)$,
where $\rho$ is an arbitrary partition in $NC(n)$:
\begin{equation}  \label{eqn:3.042}
\cf_{(i_1, \ldots , i_n); \rho} (g) 
 = \sum_{ \begin{array}{c}
{\scriptstyle  \pi \in NC(n),}  \\
{\scriptstyle \pi \leqleq \rho}
\end{array}  } \ \cf_{(i_1, \ldots , i_n); \pi} (f).
\end{equation}
\end{proposition}
\hfill $\square$

\vspace{6pt}

\begin{remark}  \label{rem:3.45}
{\rm 

$1^o$ For the proof of the above formulas (\ref{eqn:3.041}) and
(\ref{eqn:3.042}) we refer to Proposition 3.9 of \cite{BN06}.
Let us mention here that the same Proposition 3.9 of \cite{BN06}
also gives an explicit formula for how Equation (\ref{eqn:3.041})
can be inverted in order to write the coefficients of $f$ in terms 
of those of $g$. This latter fomula says that for every $n \geq 1$
and $1 \leq i_1, \ldots , i_n \leq k$ we have:
\begin{equation}  \label{eqn:3.451}
\cf_{(i_1, \ldots , i_n)} (f) 
 = \sum_{ \begin{array}{c}
{\scriptstyle  \pi \in NC(n),}  \\
{\scriptstyle \pi \leqleq 1_n}
\end{array}  } \ (-1)^{1+ | \pi |} 
\cf_{(i_1, \ldots , i_n); \pi} (g).
\end{equation}
Note that, since
$(-1)^{1+ | \pi |} \cf_{(i_1, \ldots , i_n); \pi} (g)$ can also
be written as ``$- \cf_{(i_1, \ldots , i_n); \pi} (-g)$'', an 
equivalent way of recording the formula (\ref{eqn:3.451}) is by
stating that 
\begin{equation}  \label{eqn:3.452}
\Reta^{-1} (g) = - \Reta ( -g ), \ \ \forall \, g \in \ncserk .
\end{equation}

$2^o$ Let $f,g$ be two series in $\ncserk$ such that 
$\Reta (f) = g$. An immediate consequence of Equation 
(\ref{eqn:3.041}) is that the linear and quadratic coefficients 
of $g$ are identical with the corresponding coefficients of $f$:
\[
\cf_{(i)} (g) = \cf_{(i)} (f) \mbox{ and }
\cf_{(i_1, i_2)} (g) = \cf_{(i_1, i_2)} (f), \ \ 
\forall \, 1 \leq i,i_1,i_2 \leq k.
\]
(This is because for $n \leq 2$ the only partition $\pi \in NC(n)$
which satisfies $\pi \leqleq 1_n$ is $1_n$ itself.) The first time
when we see a difference between $f$ and $g$ is when we look at 
coefficients of order 3:
\[
\cf_{(i_1,i_2,i_3)} (g) = \cf_{(i_1,i_2,i_3)} (f) 
+ \cf_{(i_1,i_3)}(f) \cdot \cf_{(i_2)} (f), \ \
\mbox{ for } 1 \leq i_1, i_2, i_3 \leq k.
\] 

$3^o$ In Section 4 of the paper we will need a formula for the 
iterations of $\Reta$, which we derive in Proposition \ref{prop:3.6}
below. The proof of this formula is based on a property of the 
partial order $\leqleq$ which was proved in Proposition 2.13 of
\cite{BN06}, and goes as follows. }
\end{remark}

\vspace{6pt}

\begin{lemma}  \label{lemma:3.5}
Let $\pi$ be a partition in $NC(n)$ such 
that $\pi \leqleq 1_n$. For every integer $p$ satisfying 
$1 \leq p \leq | \pi |$, we have that:
\[
{ } \hspace{1.5cm}
\card \Bigl\{ \rho \in NC(n) \mid \rho \geqgeq \pi \mbox{ and } 
| \rho | = p \Bigr\} 
= \left(  \begin{array}{c}
| \pi | - 1 \\
   p    - 1
\end{array}  \right) . 
\hspace{1.5cm} \square
\]
\end{lemma}

\vspace{6pt}

\begin{proposition}   \label{prop:3.6}
Let $f$ be a series in $\ncserk$ and let $s$ be in $\bR$,
$s \neq -1$. We have
\begin{equation}  \label{eqn:3.7}
\Reta \Bigl( \, s \, \Reta (f) \, \Bigr) = 
\frac{s}{1+ s} \, \Reta \Bigl( \, (1 + s)f \, \Bigr).
\end{equation}
\end{proposition}

$\ $

\noindent
{\bf Proof.} Fix $n \geq 1$ and $1 \leq i_1, \ldots , i_n \leq k$
for which we verify the equality of the coefficients of 
$z_{i_1} \cdots z_{i_n}$ for the series on the two sides of 
Equation (\ref{eqn:3.7}). We start from the left-hand side of this
equation, and compute:
\begin{align*}
\cf_{(i_1, \ldots , i_n)}
\Bigl( \, \Reta \Bigl( \, s \Reta (f) \, \Bigr) \, \Bigr)
& = \sum_{\begin{array}{c}
{\scriptstyle \rho \in NC(n)}  \\
{\scriptstyle \rho \leqleq 1_n} 
\end{array}} \ 
\cf_{(i_1, \ldots , i_n); \rho} \Bigl( \, s \Reta (f) \, \Bigr) 
\ \mbox{ (by (\ref{eqn:3.041}))}                                  \\
& = \sum_{\begin{array}{c}
{\scriptstyle \rho \in NC(n)}  \\
{\scriptstyle \rho \leqleq 1_n} 
\end{array}} \ 
s^{ | \rho |} 
\cf_{(i_1, \ldots , i_n); \rho} \Bigl( \Reta (f) \Bigr)    \\
& = \sum_{\begin{array}{c}
{\scriptstyle \rho \in NC(n)}  \\
{\scriptstyle \rho \leqleq 1_n} 
\end{array}} \ \Bigl( \,
s^{ | \rho |} \cdot
\sum_{\begin{array}{c}
{\scriptstyle \pi \in NC(n)}  \\
{\scriptstyle \pi \leqleq \rho} 
\end{array}} \ 
\cf_{(i_1, \ldots , i_n); \pi} (f) \, \Bigr) 
\ \mbox{ (by (\ref{eqn:3.042})).}                 
\end{align*}
By reversing the order of summation in the double sum that 
has appeared, we continue our sequence of equalities with:
\begin{align*}
{  }  
& = \sum_{\begin{array}{c}
{\scriptstyle \pi \in NC(n)}  \\
{\scriptstyle \pi \leqleq 1_n} 
\end{array}} \ \cf_{(i_1, \ldots , i_n); \pi} (f) \cdot
\Bigl( \,   \sum_{\begin{array}{c}
{\scriptstyle \rho \in NC(n) \ such}  \\
{\scriptstyle that \ \pi \leqleq \rho \leqleq 1_n} 
\end{array}} \ s^{ | \rho |} \, \Bigr)                          \\
& = \sum_{\begin{array}{c}
{\scriptstyle \pi \in NC(n)}  \\
{\scriptstyle \pi \leqleq 1_n} 
\end{array}} \ \cf_{(i_1, \ldots , i_n); \pi} (f) \cdot
\left(  \sum_{p=1}^{| \pi |}
{{| \pi | -1} \choose {p-1}} s^p \right) 
\ \ \mbox{ (by Lemma \ref{lemma:3.5})}                          \\
& = \sum_{\begin{array}{c}
{\scriptstyle \pi \in NC(n)}  \\
{\scriptstyle \pi \leqleq 1_n} 
\end{array}} \ s (1+ s)^{| \pi |-1} 
\cf_{(i_1, \ldots , i_n); \pi} ( f )                           \\
& = \frac{s}{1+ s} \cdot
\sum_{\begin{array}{c}
{\scriptstyle \pi \in NC(n)}  \\
{\scriptstyle \pi \leqleq 1_n} 
\end{array}} \ \cf_{(i_1, \ldots , i_n); \pi} 
\Bigl( \, (1+ s )f \, \Bigr)                                   \\
& = \cf_{(i_1, \ldots , i_n)}
\Bigl( \, \frac{s}{1+ s} \,
\Reta  \Bigl( \, (1+ s ) f \, \Bigr) 
\ \ \mbox{ (by (\ref{eqn:3.041})).}
\hspace{2cm} \square 
\end{align*}

$\ $

\begin{remark}  \label{rem:3.7}
{\rm In the case when $s=-1$, the expression 
``$\Reta \Bigl( \, s\Reta (f) \, \Bigr)$'' is not treated by the 
preceding proposition, but rather by using Equation 
(\ref{eqn:3.452}) of Remark \ref{rem:3.45}.1, which gives us that 
\begin{equation}  \label{eqn:3.071}
\Reta \Bigl( \, - \Reta (f) \, \Bigr) = - f, 
\ \ \forall \, f \in \Dalg (k).
\end{equation}  }
\end{remark}

$\ $

$\ $

\setcounter{section}{4}
\begin{center}
{\large\bf 4. The transformations \boldmath{$\bB_t$} on 
             \boldmath{$\Dalg (k)$} } 
\end{center}

\setcounter{equation}{0}
\setcounter{theorem}{0}

\begin{definition}  \label{def:4.1}
For every $t \geq 0$ define a transformation 
$\bB_t : \Dalg (k) \to \Dalg (k)$ by the formula
\begin{equation}  \label{eqn:4.1}
\bB_t ( \mu ) = \Bigl( \, \mu^{\boxplus (1+t)} \, 
\Bigr)^{\uplus (1/(1+t))}, \ \ \forall \, \mu \in \Dalg (k).
\end{equation}
\end{definition}

$\ $

Every $\bB_t$ is a bijection from $\Dalg (k)$ onto 
itself (which happens because, as noticed in Remark 
\ref{rem:2.6}, both the maps 
$\Dalg (k) \ni \mu \mapsto \mu^{\boxplus (1+t)} \in \Dalg (k)$ and 
$\Dalg (k) \ni \nu \mapsto \nu^{\uplus 1/(1+t)} \in \Dalg (k)$ are
bijective). The transformations $\{ \bB_t \mid t \geq 0 \}$ form 
in fact a semigroup under composition; this will follow from a
``commutation relation'', stated in the next proposition, satisfied 
by the convolution powers with respect to $\boxplus$ and to $\uplus$.

$\ $

\begin{proposition}  \label{prop:4.2}
Let $p,q$ be two real numbers such that $p \geq 1$ and 
$q > (p-1)/p$. We have
\begin{equation} \label{eqn:4.2}
\Bigl( \, \mu^{\boxplus p} \, \Bigr)^{\uplus q} = 
\Bigl( \, \mu^{\uplus q'} \, \Bigr)^{\boxplus p'}, \ \ 
\forall \, \mu \in {\mathcal M},
\end{equation}
where the new exponents $p', q' > 0$ are defined by
\begin{equation} \label{eqn:4.3}
p' := pq/(1-p+pq), \ \ q' := 1-p+pq.
\end{equation}
\end{proposition}

$\ $

\noindent
{\bf Proof.} If $q=1$ then it follows that $q'=1$ and $p'=p$, 
and both sides of Equation (\ref{eqn:4.2}) are equal to 
$\mu^{\boxplus p}$. For the rest of the proof we will assume that
$q \neq 1$, which implies that $q' \neq 1$ as well. Our strategy 
is to prove that the distributions on the two sides
of Equation (\ref{eqn:4.2}) have equal $R$-transforms. We prove 
this by calculating explicitly the $R$-transforms in question, where
we take advantage of the fact that the convolution powers with
respect to $\boxplus$ and with respect to $\uplus$ are scaled by
the $R$-transform and respectively by the $\eta$-series (Equations
(\ref{eqn:2.061}) and (\ref{eqn:2.063}) in Remark \ref{rem:2.6}).
The calculations may occasionally come to the point where we deal 
with the $R$-transform of a $\uplus$-power, or with the $\eta$-series
of a $\boxplus$-power; in such a situation we apply $\Reta$ (or
$\Reta^{-1}$) and go on, remaining that the compositions of 
$\Reta$'s that arise in this way are dealt with by using 
Proposition \ref{prop:3.6}. To be specific, on the 
left-hand side of (\ref{eqn:4.2}) we calculate:
\begin{align*} 
R_{ ( \mu^{\boxplus p} )^{\uplus q} } 
& = \Reta^{-1} \Bigl( \, 
    \eta_{ ( \mu^{\boxplus p} )^{\uplus q} } \, \Bigr)  
= \Reta^{-1} \Bigl( \, q \cdot \eta_{\mu^{\boxplus p}} \, \Bigr)  
\ \ \mbox{ (by (\ref{eqn:2.063})) }                               \\
& = \Reta^{-1} \Bigl( \, q \cdot 
    \Reta ( R_{\mu^{\boxplus p}} ) \, \Bigr)            
= \Reta^{-1} \Bigl( \, q \cdot \Reta ( p R_{\mu} ) \, \Bigr)     
\ \ \mbox{ (by (\ref{eqn:2.061})) }                               \\
& = - \Reta \Bigl( \, -q \cdot \Reta ( p R_{\mu} ) \, \Bigr)      
\ \ \mbox{ (by (\ref{eqn:3.452})) }                               \\
& = \frac{q}{1-q} \Reta \Bigl( \, (1-q)p \, R_{\mu} \, \Bigr) 
\ \ \mbox{ (by (\ref{eqn:3.7})). }                               
\end{align*}
On the right-hand side of (\ref{eqn:4.2}) we calculate:
\begin{align*} 
R_{ ( \mu^{\uplus q'} )^{\boxplus p'} } 
& = p' \, R_{\mu^{\uplus q'}} \hspace{2cm}
\ \ \mbox{ (by (\ref{eqn:2.061})) }                               \\
& = p' \, \Reta^{-1} \Bigl(  \eta_{\mu^{\uplus q'}} \Bigr)  
= p' \, \Reta^{-1} \Bigl(  q' \eta_{\mu} \Bigr)           
\ \ \mbox{ (by (\ref{eqn:2.063})) }                               \\
& = p' \, \Reta^{-1} \Bigl(  q' \Reta ( R_{\mu} ) \Bigr)    
= -p' \, \Reta \Bigl( -q' \Reta ( R_{\mu} ) \Bigr)       
\ \ \mbox{ (by (\ref{eqn:3.452})) }                               \\
& = (-p') \frac{-q'}{1-q'} \cdot 
    \Reta \Bigl( (1-q') \, R_{\mu} \Bigr) 
\ \ \mbox{ (by (\ref{eqn:3.7})). }                                
\end{align*}
It only remains to observe that the definition of $p'$ and $q'$ 
ensures that $1-q' = (1-q)p$ and $p'q'/(1-q') = q/(1-q)$, hence 
the two $R$-transforms calculated above are indeed equal to 
each other.
\hfill$\square$

$\ $

\begin{corollary}  \label{cor:4.3}
We have that
$\bB_s \circ \bB_t = \bB_{s+t}, \ \ 
\forall \, s,t \geq 0$.
\end{corollary}

$\ $

\noindent
{\bf Proof.} For every $s,t \geq 0$ and $\mu \in \Dalg (k)$ 
we have
\begin{eqnarray*}
\bB_s (\bB_t(\mu)) 
& = & \bB_s \left( (\mu^{\boxplus t+1})^{\uplus\frac{1}{t+1}}
            \right)                                               \\
& = & \left[ \left( (\mu^{\boxplus t+1})^{\uplus\frac{1}{t+1}}
             \right)^{\boxplus s+1} \right]^{\uplus\frac{1}{s+1}} \\
& = & \left[ \left (\mu^{\boxplus t+1})^{\boxplus
       \frac{s+t+1}{t+1}} \right)^{\uplus\frac{s+1}{s+t+1}}
             \right]^{\uplus\frac{1}{s+1}}                        \\
& = & \left( \mu^{\boxplus{s+t+1}}\right)^{\uplus
               \frac{1}{s+t+1}}                                   \\
& = & \bB_{s+t}(\mu),
\end{eqnarray*}
where at the third equality sign we used Proposition \ref{prop:4.2}
with $p = (s+t+1)/(t+1)$ and $q = (s+1)/(s+t+1)$. 
\hfill$\square$

$\ $

\begin{remark}  \label{rem:4.4}
{\rm If in the calculation for the R-transform of 
$\Bigr( \mu^{\boxplus p} \Bigl)^{\uplus q}$ that was shown 
in the proof of Proposition \ref{prop:4.2} we make $p=1+t$ and
$q = 1/(1+t)$ (for some $t>0$) we obtain 
\begin{equation}   \label{eqn:4.041}
R_{ \bB_t ( \mu ) } = \frac{1}{t} 
\Reta ( t R_{\mu} ), \ \ \forall \, \mu \in \Dalg (k),
\ \forall \, t >0.
\end{equation}
We leave it as an exercise to the reader to check that the similar
calculation done with $\eta$-series instead of $R$-transforms 
leads to the analogous formula
\begin{equation}   \label{eqn:4.042}
\eta_{ \bB_t ( \mu ) } = \frac{1}{t} 
\Reta ( t \eta_{\mu} ), \ \ \forall \, \mu \in \Dalg (k),
\ \forall \, t >0.
\end{equation}   }
\end{remark}

$\ $

\begin{remark}  \label{rem:4.5}
(Relation to the Boolean Bercovici-Pata bijection $\bB$
from \cite{BN06}.)
{\rm
In \cite{BN06} we studied a bijection 
$\bB : \Dalg (k) \to \Dalg (k)$ defined via the requirement that
\begin{equation}
R_{\bB ( \mu )} = \eta_{\mu}, \ \ \forall \, \mu \in \Dalg (k).
\end{equation}
It is immediate that $\bB$ coincides with the transformation 
$\bB_1$ obtained by putting $t=1$ in Definition \ref{def:4.1}.
Indeed, for every $\mu \in \Dalg (k)$ we have
\begin{align*}
R_{\bB_1 ( \mu )}
& = \Reta ( R_{\mu} ) \ \ \mbox{ (by making $t=1$ in 
                                  Equation (\ref{eqn:4.041}))} \\
& = \eta_{\mu} \ \ \mbox{ (by definition of $\Reta$); }     
\end{align*}
this implies that $\bB_1 ( \mu ) = \bB ( \mu )$, since 
$\bB_1 ( \mu )$ and $\bB ( \mu )$ have the same $R$-transform.  }
\end{remark}

\vspace{6pt}

\begin{remark}  \label{rem:4.6}
{\rm 
An intriguing property of the map $\bB$ which was observed in 
\cite{BN06} is that it is a homomorphism with respect to the 
operation of {\em free multiplicative} convolution $\boxtimes$
on $\Dalg (k)$. This operation is defined as follows. Given 
$\mu, \nu \in \Dalg (k)$, one can always find random variables 
$x_1, \ldots , x_k, y_1, \ldots , y_k$ in a non-commutative 
probability space $( \cA, \varphi )$ such that the joint 
distribution of the $k$-tuple $x_1, \ldots , x_k$ is equal 
to $\mu$, the joint distribution of the 
$k$-tuple $y_1, \ldots , y_k$ is equal to $\nu$, and such that 
$\{ x_1, \ldots , x_k \}$ is freely independent from 
$\{ y_1, \ldots , y_k \}$ in $( \cA , \varphi )$. The joint 
distribution of the $k$-tuple $x_1y_1, \ldots , x_ky_k$ turns 
out to depend only 
on $\mu$ and $\nu$; and the free multiplicative convolution 
$\mu \freetimes \nu$ is equal, by definition, to the joint 
distribution of $x_1y_1, \ldots , x_ky_k$. 

In the remaining part of this section we will show that every 
$\bB_t$ is a homomorphism with respect to $\boxtimes$. The argument 
is short, because it takes advantage of what had already been 
proved in \cite{BN06} -- the essential point is to use Theorem 7.3 
of that paper. We mention that in the 1-dimensional case another 
derivation of the $\boxtimes$-homomorphism property of $\bB_t$
can be obtained by using the concept of S-transform (see Section 3 
of \cite{BN07}). 

In the proof that $\bB_t$ is a $\boxtimes$-homomorphism we will
also use a binary operation denoted by $\freestar$ on 
$\ncserk$, which was introduced in \cite{NS96}, and is uniquely 
determined by the fact that 
\begin{equation}  \label{eqn:4.061}
R_{\mu} \ \freestar \ R_{\nu} = R_{\mu \boxtimes \nu}, 
\ \ \forall \, \mu , \nu \in \Dalg (k).
\end{equation}
In other words, $\freestar$ is the operation with formal power 
series which reflects the multiplication of two free $k$-tuples
in terms of their $R$-transforms.

A remarkable fact proved in 
Theorem 7.3 of \cite{BN06} is that we also have
\begin{equation} \label{eqn:4.062}
\eta_{\mu} \ \freestar \ \eta_{\nu} = \eta_{\mu \boxtimes \nu}, 
\ \ \forall \, \mu , \nu \in \Dalg (k).
\end{equation}
That is, $\freestar$ is at the same time the operation with formal 
power series which reflects the multiplication of two free 
$k$-tuples in terms of their $\eta$-series. It is immediate
that formula (\ref{eqn:4.062}) is actually just another form 
of stating the $\boxtimes$-multiplicativity of $\bB$. We prefer 
this formula which makes explicit use of $\freestar$, 
because we want to combine it with other properties that 
$\freestar$ has, in connection to dilations and scalar 
multiplication of power series (as reviewed in Remark \ref{rem:4.8}
below).  }
\end{remark}

\vspace{6pt}

\begin{definition}  \label{def:4.7}
{\rm
$1^o$ For $\mu \in \Dalg (k)$ and $r>0$ we denote by $\mu \circ D_r$
the distribution in $\Dalg (k)$ determined by the condition that
\begin{equation}
( \mu \circ D_r ) (X_{i_1} \cdots X_{i_n})
= r^n \cdot \mu (X_{i_1} \cdots X_{i_n}), \ \ 
\forall \, n \geq 1, \ \forall \, 
1 \leq i_1, \ldots , i_n \leq k.
\end{equation}
$\mu \circ D_r$ is called the {\em dilation of $\mu$ by $r$}.

$2^o$ For $f \in \ncserk$ and $r>0$ we denote by $f \circ D_r$
the series in $\ncserk$ determined by the condition that
\begin{equation}
\cf_{(i_1, \ldots , i_n)} (f \circ D_r) =
r^n \cdot \cf_{(i_1, \ldots , i_n)} (f),
\ \ \forall \, n \geq 1, \ \forall \, 
1 \leq i_1, \ldots , i_n \leq k.
\end{equation}                           
$f \circ D_r$ is called the {\em dilation of $f$ by $r$}.  }
\end{definition}

\vspace{6pt}

\begin{remark}  \label{rem:4.8}
{\rm It is easy to see, directly from the definitions, that all 
three series $M_{\mu}$, $R_{\mu}$, $\eta_{\mu}$ associated to a 
distribution $\mu \in \Dalg (k)$ behave well with respect to 
dilations; that is, we have
\begin{equation}  \label{eqn:4.081}
M_{\mu \circ D_r} = M_{\mu} \circ D_r, \ \
R_{\mu \circ D_r} = R_{\mu} \circ D_r, \ \
\eta_{\mu \circ D_r} = \eta_{\mu} \circ D_r, \ \ 
\forall \, \mu \in \Dalg (k), \ \forall \, r > 0.
\end{equation}
Let us also record here two formulas from \cite{NS96} which 
involve dilations and the operation $\freestar$. The first formula
simply says that $\freestar$ behaves well with respect to 
dilations:
\begin{equation}  \label{eqn:4.082}
(f \circ D_r) \ \freestar \ g = 
f \ \freestar \ (g \circ D_r) = 
(f \ \freestar \ g) \circ D_r, \ \ 
\forall \, f,g \in \ncserk , \ \forall \, r>0. 
\end{equation}
The second formula puts into evidence a special connection with
scalar multiplication of series. While $\freestar$ is highly 
non-linear (and doesn't generally behave well with respect to 
scalar multiplication), it is remarkable that we have
\begin{equation}  \label{eqn:4.083}
( rf) \ \freestar \ ( rg) = 
r \Bigl( \, (f \ \freestar \ g) \circ D_r \, \Bigr),
\forall \, f,g \in \ncserk , \ \forall \, r>0. 
\end{equation}
For the proof of (\ref{eqn:4.082}) and (\ref{eqn:4.083}) we refer 
to Notation 4.1 and Lemma 4.4 of \cite{NS96}.   }
\end{remark}

$\ $

In order to prove that $\bB_t$ is a homomorphism with respect 
to $\boxtimes$ we will show in the next proposition that, in fact,
each of the two kinds of convolution powers involved in the 
definition of $\bB_t$ is ``only a dilation away'' from being 
itself a $\boxtimes$-homomorphism.

$\ $

\begin{proposition}  \label{prop:4.9}
For every $t>0$ and every $\mu , \nu \in \Dalg (k)$ we have
\begin{equation}  \label{eqn:4.091}
( \mu^{\boxplus t} ) \boxtimes ( \nu^{\boxplus t} )  = 
( \mu \boxtimes \nu )^{\boxplus t} \circ D_t
\end{equation}
and
\begin{equation}  \label{eqn:4.092}
( \mu^{\uplus t} ) \boxtimes ( \nu^{\uplus t} ) = 
( \mu \boxtimes \nu )^{\uplus t} \circ D_t.
\end{equation}       
\end{proposition}

$\ $

\noindent
{\bf Proof.}
In order to establish the fomula (\ref{eqn:4.091})
we check that the distributions appearing on the two sides of 
this formula have the same $R$-transform:
\begin{align*}
R_{( \mu^{\boxplus t} ) \boxtimes ( \nu^{\boxplus t} ) }
& = R_{\mu^{\boxplus t}} \ \freestar \ R_{\nu^{\boxplus t}} 
\ \ \mbox{ (by (\ref{eqn:4.061})) }                         \\
& = (tR_{\mu}) \ \freestar \ (tR_{\nu})      
\ \ \mbox{ (by (\ref{eqn:2.061})) }                         \\
& = \Bigl( t \, (R_{\mu} \ \freestar \ R_{\nu}) \Bigr)
\circ D_t  \ \ \mbox{ (by (\ref{eqn:4.083})) }              \\
& = \Bigl( t \, R_{\mu \boxtimes \nu} \Bigr) \circ D_t      
\ \ \mbox{ (by (\ref{eqn:4.061})) }                         \\
& = \Bigl( R_{( \mu \boxtimes \nu )^{\boxplus t}} \Bigr)
\circ D_t  \ \ \mbox{ (by (\ref{eqn:2.061})) }              \\
& = R_{( \mu \boxtimes \nu )^{\boxplus t} \circ D_t} 
\ \ \mbox{ (by (\ref{eqn:4.081})). } 
\end{align*}

The verification of (\ref{eqn:4.092}) is done in 
a similar way, where now we check that the 
distributions on the two sides of the formula have
identical $\eta$-series. The calculation is virtually
identical to the one shown in the verification of 
(\ref{eqn:4.091}), only that we have to replace everywhere 
$R$-transforms by $\eta$-series, and $\boxplus$-powers 
by $\uplus$-powers. (An important point included in this 
``mutatis mutandis'' argument is that, right at the beginning 
of the calculation, we can invoke the formula (\ref{eqn:4.062}) 
relating $\eta$-series to the operation $\freestar$.)
\hfill$\square$

$\ $

\begin{corollary}  \label{cor:4.10}
For every $t \geq 0$, the transformation $\bB_t$ of $\Dalg (k)$ 
is a homomorphism for $\boxtimes$. That is, we have 
\begin{equation}  \label{eqn:4.101}
\bB_t ( \mu \freetimes \nu ) = \bB_t ( \mu ) 
\boxtimes \bB_t ( \nu ), \ \ \forall \, \mu , \nu \in \Dalg (k).
\end{equation}  
\end{corollary}

$\ $

\noindent
{\bf Proof.} This is a straightforward 
consequence of Proposition \ref{prop:4.9}: the dilation factors 
which appear when we take succesively the powers ``$\boxplus (t+1)$'' 
and ``$\uplus 1/(t+1)$'' cancel each other, and we are left with the 
plain $\boxtimes$-multiplicativity stated in Equation 
(\ref{eqn:4.101}). 
\hfill$\square$

$\ $

The results of this section are thus summarized in the 
following theorem, which puts together Corollary \ref{cor:4.3},
Remark \ref{rem:4.5}, and Corollary \ref{cor:4.10}.

$\ $

\begin{theorem}   \label{thm:4.11}
The bijections $\bB_t : \Dalg (k) \to \Dalg (k)$ introduced 
in Definition \ref{def:4.1} have the following properties:

$1^o$ $\bB_s \circ \bB_t = \bB_{s+t}$, for every $s,t \geq 0$.

$2^o$ $\bB_1 = \bB$, the multi-variable Boolean Bercovici-Pata
bijection introduced in \cite{BN06}.

$3^o$ Every $\bB_t$ is a homomorphism for the free multiplicative
convolution $\boxtimes$ on $\Dalg (k)$.
\end{theorem}
\hfill $\square$

$\ $

$\ $

\setcounter{section}{5}
\begin{center}
{\large\bf 5. A formula for the moments of the free Brownian motion}
\end{center}

\setcounter{equation}{0}
\setcounter{theorem}{0}
Our goal in this section is to prove an explicit formula via
summations over non-crossing partitions for moments 
$( \nu \boxplus \gamma_t ) ( X_{i_1} \cdots X_{i_n} )$, 
where $\nu$ is an arbitrary distribution in $\Dalg (k)$ and 
$\gamma_t$ is defined as follows.

$\ $

\begin{notation}  \label{def:5.1}
{\rm 
For $t > 0$ we will denote by $\gamma_t \in \Dalg (k)$ 
the joint distribution of a $k$-tuple $(x_1, \ldots , x_k)$ where
$x_1, \ldots , x_k$ form a free family, and every $x_i$ has a 
centered semicircular distribution of variance $t$.   }
\end{notation}

$\ $

The formula for moments which is the main result of the section 
will be stated in Proposition \ref{prop:5.5}. We start by 
introducing a few natural conventions of notations for 
non-crossing partitions that will be useful in Proposition
\ref{prop:5.5}.

$\ $

\begin{remark}  \label{rem:5.2}

{\rm 
$1^o$ It will be convenient that instead of sticking strictly
to ``$NC(n)$'', we use the more general notation ``$NC(M)$''
for an arbitrary totally ordered finite set $M$. Of course,
$NC(M)$ can always be identified canonically to $NC( \, |M| \, )$, 
where one uses the unique increasing bijection from $M$ onto 
$\{ 1, \ldots , |M| \}$ in order to identify partitions of $M$ 
with partitions of $\{ 1, \ldots , |M| \}$.

$2^o$ Let $M$ be a totally ordered finite set, and let $L$ be a 
non-empty subset of $M$. For $\pi \in NC(M)$ we can consider the 
restricted partition $\pi \mid L$ of $L$ into blocks of the 
form $A \cap L$, with $A$ block of $\pi$ such that 
$A \cap L \neq \emptyset$. It is immediately verified 
that $\pi \mid L \in NC(L)$ (where 
$L$ is endowed with the total order inherited from $M$).

$3^o$ Let $M$ be a totally ordered finite set, and suppose that 
$M = L_1 \cup L_2$, disjoint union. If $\pi_1$ is a partition of 
$L_1$ and $\pi_2$ is a partition of $L_2$, 
then there is an obvious way
of putting $\pi_1$ and $\pi_2$ together to form a partition of
$M$; we will denote this partition by $\pi_1 \sqcup \pi_2$. It is 
clear that in order to have $\pi_1 \sqcup \pi_2 \in NC(M)$ it is 
necessary but not sufficient that $\pi_1 \in NC(L_1)$ and 
$\pi_2 \in NC(L_2)$.

$4^o$ Let $M, L_1, L_2$ be as above and let $\pi_1$ be a fixed 
partition in $NC(L_1)$. It is easy to see that among the 
partitions $\pi_2 \in NC(L_2)$ with the property that 
$\pi_1 \sqcup \pi_2 \in NC(M)$ there is one, $\widehat{\pi}$,
which is larger than all the others with respect to reversed 
refinement order on $NC(L_2)$. So $\widehat{\pi} \in NC( L_2 )$
is characterized by the fact that for a partition 
$\pi_2 \in NC(L_2)$ we have the equivalence
\begin{equation}  \label{eqn:5.021}
\pi_1 \sqcup \pi_2 \in NC(A) \ \Leftrightarrow \
\pi_2 \leq \widehat{\pi}.
\end{equation}                  }
\end{remark}

$\ $

The formula for moments that will be proved in Proposition 
\ref{prop:5.5} uses the class of non-crossing partitions 
discussed in the following notation.

$\ $

\begin{notation}  \label{def:5.3}
{\rm  Let $n$ be a positive integer.

$1^o$ We will denote by $NC_{\leq 2} (n)$ the set of partitions 
$\rho \in NC(n)$ such that every block of $\rho$ has either 1 or 
2 elements.         

$2^o$ For a partition $\rho$ in $NC_{\leq 2}(n)$ we will denote 
by $D( \rho )$ the union of all doubletons (2-element blocks) of 
$\rho$, and by $S( \rho )$ the union of all singletons (1-element 
blocks) of $\rho$. Thus $D( \rho ) \cup S( \rho )$ = 
$\{ 1, \ldots , n \}$ (disjoint union). 

$3^o$ Let $\rho$ be in $NC_{\leq 2}(n)$, and let us consider 
the partition $\rho \mid D( \rho ) \in NC( \, D( \rho ) \, )$. 
We will denote by $\widehat{\rho}$ the maximal partition in 
$NC( \, S( \rho ) \, )$ that can be combined with 
$\rho \mid D( \rho )$ into a non-crossing partition of 
$\{ 1, \ldots , n \}$, in the sense discussed in part $4^o$ of
the preceding remark. Thus $\widehat{\rho}$ is characterized
by the fact that for a partition 
$\sigma \in NC( \, S( \rho ) \, )$ we have the equivalence
\begin{equation}  \label{eqn:5.031}
\Bigl( \, \rho \mid D( \rho ) \, \Bigr) \sqcup
\sigma \in NC(n) \ \Leftrightarrow \
\sigma \leq \widehat{\rho}.
\end{equation}                  

\vspace{6pt}

\noindent
[A concrete example illustrating the parts $2^o$ and $3^o$ of this 
notation: say that $n=9$ and that 
\begin{equation}  \label{eqn:5.25}
\rho = \Bigl\{ \, \{ 1 \}, \ \{ 2,8 \} , \{ 3 \} , 
\{ 4,5 \} , \{ 6 \} , \{ 7 \} , \{ 9 \} \, \Bigr\}
\in NC_{\leq 2} (9).
\end{equation}
Then $D( \rho ) = \{ 2,4,5,8 \}$,
$S( \rho ) = \{ 1,3,6,7,9 \}$, and we have 
$\rho \mid D( \rho ) = \{ \, \{ 2,8 \}, \, \{ 4,5 \} \, \} 
\in NC( \, D( \rho ) \, )$ and 
$\widehat{\rho} = \{ \ \{ 1,9 \}, \ \{ 3,6,7 \} \, \} \in
NC( \, S( \rho ) \, )$.]    }
\end{notation}

\vspace{6pt}

\begin{proposition}  \label{prop:5.5}
Let $\nu$ be a distribution in $\Dalg (k)$, and let $\gamma_t$ be
as described in Notation \ref{def:5.1}. For every $n \geq 1$ and
$1 \leq i_1, \ldots , i_n \leq k$ we have
\begin{equation}   \label{eqn:5.051}
( \nu \boxplus \gamma_t ) ( X_{i_1} \cdots X_{i_n} ) = 
\hspace{5cm}  {  }
\end{equation}
\[
\sum_{\rho \in NC_{\leq 2}(n)} \
\Bigl( \, \Bigl( \, \prod_{  \begin{array}{c}
{\scriptstyle B \ 2-element \ block}  \\
{\scriptstyle of \ \rho, \ B=\{ p,q \} }
\end{array} } \ t \delta_{i_p, i_q} \, \Bigr) \cdot
\cf_{( \, (i_1, \ldots , i_n) | S( \rho ) \, ); \widehat{\rho}}
\ (M_{\nu}) \, \Bigr) .
\]
\end{proposition}

\vspace{6pt}

\begin{remark}  \label{rem:5.6}
{\rm Let us comment a bit on what is achieved by the formula 
(\ref{eqn:5.051}). An important point is, of course, that we
explicitly identify a combinatorial structure -- namely 
$NC_{\leq 2} (n)$ -- which indexes the sum leading to 
$( \nu \boxplus \gamma_t ) ( X_{i_1} \cdots X_{i_n} )$.
Let us moreover fix a partition $\rho \in NC_{\leq 2}(n)$ and
let us examine the term indexed by $\rho$ on the right-hand side
of (\ref{eqn:5.051}). 

First there is an issue of compatibility.
Let us say that ``$\rho$ is compatible with the $n$-tuple 
$(i_1, \ldots , i_n)$'' when the following happens: whenever
$B = \{ p,q \}$ is a 2-element block of $\rho$, it follows that
$i_p = i_q$. If $\rho$ is not compatible with 
$(i_1, \ldots , i_n)$, then the term indexed by $\rho$ on the 
right-hand side of (\ref{eqn:5.051}) vanishes. 

Suppose then that $\rho$ is compatible 
with $(i_1, \ldots , i_n)$. Let 
$S( \rho ) = \{ b_1 < b_2 < \cdots < b_m \}$ be the
set of singletons of $\rho$, and let $\widehat{\rho}$ be the 
non-crossing partition of $S( \rho )$ that was put into evidence 
in Notation \ref{def:5.3}.3. The term indexed by $\rho$ on the 
right-hand side of (\ref{eqn:5.051}) is then equal to
\begin{equation}  \label{eqn:5.061}
t^d \, \cf_{( i_{b_1}, i_{b_2}, \ldots , 
              i_{b_m} ); \widehat{\rho} } \ ( M_{\nu} ),
\end{equation}
where $d = (n-m)/2$ is the number of doubletons of $\rho$,
and where the generalized coefficient 
$\cf_{( i_{b_1}, i_{b_2}, \ldots , i_{b_m} ); \widehat{\rho} }$
is as in the above Definition \ref{def:3.3}. (Note the detail 
that in (\ref{eqn:5.061}) the partition $\widehat{\rho}$ is 
viewed, in the canonical way, as a partition from $NC(m)$.)

A concrete example: look again at the example of $\rho \in NC(9)$
given for illustration at the end of Notation \ref{def:5.3}. 
There we had $S( \rho ) = \{ 1,3,6,7,9 \}$, and 
$\widehat{\rho} = \{ \, \{ 1,9 \}, \ \{ 3,6,7 \} \, \}$
$\in NC( \, S( \rho ) \, )$. Thus the generalized coefficient 
of $M_{\nu}$ we have to look at is 
$\cf_{ (i_1, i_3, i_6, i_7, i_9); \widehat{\rho} } \ (M_{\nu})$,
which is just $\nu ( X_{i_1} X_{i_9} ) 
\nu ( X_{i_3} X_{i_6} X_{i_7} )$. Hence the term indexed by $\rho$
in the sum on the right-hand side of (\ref{eqn:5.051}) is in this
concrete example equal to
\[
\left\{  \begin{array}{ll}
t^2 \,  \nu ( X_{i_1} X_{i_9} ) \nu ( X_{i_3} X_{i_6} X_{i_7} )
        & \mbox{ if $i_2 = i_8$ and $i_4 = i_5$ }               \\
0       & \mbox{ otherwise. }
\end{array}  \right.
\]                                    }
\end{remark}

\vspace{10pt}

\begin{remark}  \label{rem:5.7}
{\rm 
We now move towards proving the formula stated in 
Proposition \ref{prop:5.5}. In preparation of the proof, let us 
review the basic ``moments vs. free cumulants'' formula which 
expresses the moments of a distribution $\mu \in \Dalg (k)$ 
in terms of its free cumulants -- that is, in terms of the 
coefficients of the $R$-transform $R_{\mu}$.
This formula says that
\begin{equation}   \label{eqn:5.071} 
\cf_{(i_1, \ldots , i_n)} (M_{\mu}) = \sum_{\pi \in NC(n)} 
\cf_{(i_1, \ldots , i_n); \pi} (R_{\mu}), 
\ \ \forall \, n \geq 1, \ \forall \,
1 \leq i_1, \ldots , i_n \leq k;
\end{equation}
and more generally, that for any $\rho \in NC(n)$ we have
\begin{equation}   \label{eqn:5.072} 
\cf_{(i_1, \ldots , i_n); \rho} ( M_{\mu} )
 = \sum_{ \begin{array}{c}
{\scriptstyle \pi \in NC(n),}  \\
{\scriptstyle \pi \leq \rho} 
\end{array}  } \
\cf_{(i_1, \ldots , i_n); \pi} ( R_{\mu} )
\end{equation}
(where Equation (\ref{eqn:5.071}) corresponds to the case when 
$\rho = 1_n$). For more details on this, see Lectures 11 and 16 
of \cite{NS06}. 

Also in preparation of the proof of Proposition \ref{prop:5.5}
it is convenient to introduce the following elements of 
notation.  }
\end{remark}

\vspace{10pt}

\begin{notation}  \label{def:5.8}
{\rm Let $n$ be a positive integer.

$1^o$ For $\rho \in NC_{\leq 2}(n)$ and $\pi \in NC(n)$ we 
will write ``$ \rho \tleft \pi$'' to mean that every
2-element block $B$ of $\rho$ also is a block of $\pi$.
(That is: if the 2-element blocks of $\rho$ are
$B_1, \ldots , B_p$, then $\pi$ must be of the form 
$\pi = \{ B_1, \ldots , B_p, C_1, \ldots , C_q \}$, with
$q \geq 0$ and $C_1 \cup \cdots \cup C_q = S( \rho )$.)

$2^o$ Let $i_1, \ldots , i_n$ be some indices in 
$\{ 1, \ldots , k \}$. We will denote by 
$NC_{\leq 2}(n; i_1, \ldots , i_n)$ the set of partitions 
$\rho \in NC_{\leq 2} (n)$ with the property that whenever 
$B = \{ p,q \}$ is a 2-element block of $\rho$, it follows that 
$i_p = i_q$.  }
\end{notation}  

$\ $

\noindent
{\bf Proof of Proposition \ref{prop:5.5}.} 
We fix for the whole proof a positive integer $n$ and 
some indices $1 \leq i_1, \ldots , i_n \leq k$ for which 
we will prove that Equation (\ref{eqn:5.051}) holds. 

We start from the left-hand side of the equation. From the 
moment-cumulant formula (\ref{eqn:5.071}) and the fact that 
$R_{\nu \boxplus \gamma_t} = R_{\nu} + R_{\gamma_t}$, we 
have:
\begin{equation}  \label{eqn:5.7}
( \nu \boxplus \gamma_t ) (X_{i_1} \cdots X_{i_n})
= \sum_{\pi \in NC(n)} 
\cf_{(i_1, \ldots , i_n ); \pi} ( R_{\nu} + R_{\gamma_t} ).
\end{equation}

Now let us fix for the moment a partition $\pi \in NC(n)$, and 
let us look at the term indexed by $\pi$ on the right-hand side
of (\ref{eqn:5.7}). We write this term explicitly:
\begin{equation}  \label{eqn:5.8}
\cf_{(i_1, \ldots , i_n ); \pi} ( R_{\nu} + R_{\gamma_t} )
= \prod_{ \begin{array}{c}
{\scriptstyle A \ block}  \\
{\scriptstyle of \ \pi}
\end{array}  }   \Bigr( \,
\cf_{(i_1, \ldots , i_n) \vert A} (R_{\nu}) +
\cf_{(i_1, \ldots , i_n) \vert A} (R_{\gamma_t} ) \, \Bigr)
\end{equation}
and we expand the product on the right-hand side of 
(\ref{eqn:5.8}) into a sum of $2^{| \pi |}$ terms. The general 
term of the sum is obtained by splitting the set of blocks of 
$\pi$ into two sets
of blocks $\cS_1$ and $\cS_2$, and by forming the product
\begin{equation}  \label{eqn:5.9}
\Bigl( \, \prod_{ A \in \cS_1}
\cf_{(i_1, \ldots , i_n) \vert A} (R_{\nu}) \, \Bigr) \cdot
\Bigl( \, \prod_{ B \in \cS_2}
\cf_{(i_1, \ldots , i_n) \vert B} (R_{\gamma_t}) \, \Bigr) .
\end{equation}
But a fundamental fact about free semicircular systems is that 
the $R$-transform of $\gamma_t$ is just
\[
R_{\gamma_t} (z_1, \ldots , z_k) = t(z_1^2 + \cdots + z_k^2)
\]
(see \cite{NS06}, Lectures 11 and 16).
Thus the second product in (\ref{eqn:5.9}) is non-zero if and
only if every block $B \in \cS_2$ is of the form $B= \{ p,q \}$
with $1 \leq p<q \leq n$ such that $i_p = i_q$. When this 
requirement is satisfied, the set $\cS_2$ of blocks of $\pi$ 
corresponds naturally to a partial pairing $\rho \in 
NC_{\leq 2} (n; i_1, \ldots , i_n)$ such that $\rho \tleft \pi$
(where Notation \ref{def:5.8} is used). For our fixed 
$\pi \in NC(n)$ we thus arrive to an equation of the form
\begin{equation}  \label{eqn:5.10}
\cf_{(i_1, \ldots , i_n ); \pi} ( R_{\nu} + R_{\gamma_t} )
= \sum_{\begin{array}{c}
{\scriptstyle \rho \in NC_{\leq 2} (n;i_1, \ldots , i_n)}  \\
{\scriptstyle such \ that \ \rho \tleft \pi}
\end{array} } \ \term_{\rho}
\end{equation}
where the quantities ``$\term_{\rho}$'' are further discussed 
in the next paragraph.

So let $\pi \in NC(n)$ be as in the preceding paragraph, and let
$\rho \in NC_{\leq 2} (n; i_1, \ldots , i_n)$ be such 
that $\rho \tleft \pi$. In connection to this $\rho$ we will use
the notations $D(\rho), S( \rho )$ and 
$\widehat{\rho} \in NC( \, S( \rho ) \, )$ that were introduced 
in Notation \ref{def:5.3}. The contribution ``$\term_{\rho}$'' to 
the sum (\ref{eqn:5.10}) is of the form shown in (\ref{eqn:5.9}),
where $\cS_2$ is the set of blocks of $\pi$ which also are 
2-element blocks of $\rho$. The product 
``$\prod_{B \in \cS_2} \cdots$'' in (\ref{eqn:5.9}) is then 
clearly equal to $t^{|D( \rho )|/2}$. For the other product 
``$\prod_{A \in \cS_1} \cdots$'' in (\ref{eqn:5.9}) 
we note that the union of the blocks counted
in $\cS_1$ is equal to $S( \rho )$, and this gives us that 
\[
\prod_{ A \in \cS_1}
\cf_{(i_1, \ldots , i_n) \vert A} (R_{\nu}) =
\cf_{( \, (i_1, \ldots , i_n) | S( \rho ) \, );
( \pi | S( \rho ))} \ ( R_{\nu} ).
\]

The conclusion of the preceding two paragraphs of the proof is 
that for every $\pi \in NC(n)$ we have
\begin{equation}  \label{eqn:5.11}
\cf_{(i_1, \ldots , i_n ); \pi} ( R_{\nu} + R_{\gamma_t} )
= \sum_{\begin{array}{c}
{\scriptstyle \rho \in NC_{\leq 2} (n;i_1, \ldots , i_n)}  \\
{\scriptstyle such \ that \ \rho \tleft \pi}
\end{array} } \ t^{|D( \rho )|/2} \cdot
\cf_{( \, (i_1, \ldots , i_n) | S( \rho ) \, );
( \pi | S( \rho ))} \ ( R_{\nu} ).
\end{equation}
We now sum over $\pi$ in Equation (\ref{eqn:5.11}). On the 
left-hand side the sum over $\pi$ gives us 
$( \nu \boxplus \gamma_t ) (X_{i_1} \cdots X_{i_n})$,
as we knew since (\ref{eqn:5.7}). On the right-hand side of 
(\ref{eqn:5.11}) we get a double sum, over $\pi$ and $\rho$; 
we interchange the order of summation in this double sum,
to obtain:
\begin{equation}  \label{eqn:5.12}
\sum_{\rho \in NC_{\leq 2}(n; i_1, \ldots , i_n)} \ 
t^{|D( \rho )|/2} \Bigl( \
\sum_{\begin{array}{c}
{\scriptstyle  \pi \in NC(n) \ such }  \\
{\scriptstyle that \ \rho \tleft \pi} 
\end{array} } \ 
\cf_{( \, (i_1, \ldots , i_n) | S( \rho ) \, );
( \pi | S( \rho ))} \ ( R_{\nu} ) \, \Bigr) .
\end{equation}

It is now the turn of $\rho$ to be fixed, while
we examine the summation over $\pi$ that has appeared 
in (\ref{eqn:5.12}). By taking into account the 
discussion from Notation \ref{def:5.3}, it is immediate that 
every partition $\pi \in NC(n)$ with the property that 
$\rho \tleft \pi$ is obtained in a unique way as 
$( \rho \mid D( \rho ) ) \sqcup \sigma$, where 
$\sigma \in NC( \, S( \rho ) \, )$ is such that 
$\sigma \leq \widehat{\rho}$ (see the equivalence 
(\ref{eqn:5.031}) in Notation \ref{def:5.3}). It follows that 
the inside sum over $\pi$ in (\ref{eqn:5.12}) is equal to
\[
\sum_{\begin{array}{c}
{\scriptstyle  \sigma \in NC( \, S( \rho ) \, ) \ such }  \\
{\scriptstyle that \ \sigma \leq \widehat{\rho} } 
\end{array} } \ 
\cf_{( \, (i_1, \ldots , i_n) | S( \rho ) \, ); \sigma }
\ ( R_{\nu} ).
\]
But the latter quantity is in turn equal to
$\cf_{( \, (i_1, \ldots , i_n) | S( \rho )); 
\widehat{\rho} } \ (M_{\nu})$,
due to the moments vs. free cumulant formula (used now in the 
more general form that was reviewed in (\ref{eqn:5.072})). 
Replacing this in (\ref{eqn:5.12}) takes us precisely to the 
right-hand side of Equation (\ref{eqn:5.051}), and this 
concludes the proof.
\hfill$\square$

$\ $

$\ $

\setcounter{section}{6}
\begin{center}
{\large\bf 6. Relation between \boldmath{$\bB_t$}
              and the free Brownian motion}
\end{center}

\setcounter{equation}{0}
\setcounter{theorem}{0}
Recall from Remark \ref{rem:2.5} that the map 
$\mu \mapsto \eta_{\mu}$ is a bijection from $\Dalg (k)$ 
onto the space of series $\ncserk$. It thus makes sense to define 
a map $\Phi : \Dalg (k) \to \Dalg (k)$ via the $\eta$-series
prescription described as follows.

\begin{definition}  \label{def:6.1}
{\rm For every $\nu \in \Dalg (k)$, we let $\Phi ( \nu )$ be
the unique distribution $\mu \in \Dalg (k)$ which has 
$\eta$-series given by:
\begin{equation}  \label{eqn:6.1}
\eta_{\mu} (z_1, \ldots , z_k) = \sum_{i=1}^k
z_i \Bigl( 1+ M_{\nu} (z_1, \ldots , z_k) \Bigr) z_i.
\end{equation}  }
\end{definition}

Our goal in the present section is to prove the following result.

$\ $

\begin{theorem}  \label{thm:6.2}
Let $\nu$ be a distribution in $\Dalg (k)$. We have that
\begin{equation}   \label{eqn:6.2}
\Phi ( \, \nu \boxplus \gamma_t \, ) = 
\bB_t ( \, \Phi ( \nu ) \, ), \ \ \forall \, t > 0,
\end{equation}
where $\gamma_t \in \Dalg (k)$ is the distribution of the 
scaled free semicircular system from Notation \ref{def:5.1}.
\end{theorem}

$\ $

A key point in the proof of Theorem \ref{thm:6.2} will be to use
a natural combinatorial construction of ``assigning singletons 
to doubletons'' in a partial pairing, which is described next.

\begin{remark}  \label{rem:6.3} 
(``Assigning singletons to 
doubletons for $\rho \in NC_{\leq 2}(n)$''.)
{\rm
Let a partition $\rho \in NC_{\leq 2}(n)$ be given. We will 
denote by $\alpha ( \rho )$ the non-crossing partition of 
$\{ 0, 1, \ldots , n , n+1 \}$ which is obtained as follows.
Start with the partial pairing of
$\{ 0, 1, \ldots , n , n+1 \}$ that is obtained by adding to
$\rho$ the 2-element block $\{ 0, n+1 \}$. Consider the picture
of this new partial pairing (drawn in the usual way -- with 
the points $0,1, \ldots , n , n+1$ on a horizontal line, and
with a family of non-intersecting ``hooks'' drawn under that 
horizontal line, to represent the 2-element blocks of the 
partial pairing). In this picture we draw some additional 
vertical line segments, starting at every singleton of $\rho$,
and going down until they meet a hook representing a doubleton.
When these new vertical segments are added to the picture, we 
now have the picture of a non-crossing partition of 
$\{ 0, 1, \ldots , n , n+1 \}$, which will be denoted by 
$\alpha ( \rho )$.

A concrete example: if $n=9$ and $\rho \in NC_{\leq 2}(n)$ is as 
in (\ref{eqn:5.25}) from Notation \ref{def:5.3}, then 
$\alpha ( \rho ) = \{ \ \{ 0,1,9,10 \} , \ \{ 2,3,6,7,8 \} , 
\{ 4,5 \} \ \}$, and
the pictures of $\rho$ and of $\alpha ( \rho )$ look as follows:
\[
\rho \ = \ 
\setlength{\unitlength}{0.3cm}
\begin{picture}(9,4)\thicklines
\put(0,0){\line(0,1){1}} 
\put(1,-1){\line(0,1){2}} \put(1,-1){\line(1,0){6}} 
\put(2,0){\line(0,1){1}} 
\put(3,0){\line(0,1){1}} \put(3,0){\line(1,0){1}} 
\put(4,0){\line(0,1){1}} 
\put(5,0){\line(0,1){1}} 
\put(6,0){\line(0,1){1}} 
\put(7,-1){\line(0,1){2}} 
\put(8,0){\line(0,1){1}} 
\put(-0.3,1.7){1} \put(0.7,1.7){2} \put(1.7,1.7){3} 
\put(2.7,1.7){4}  \put(3.7,1.7){5} \put(4.7,1.7){6} 
\put(5.7,1.7){7}  \put(6.7,1.7){8} \put(7.7,1.7){9}
\end{picture} 
\ \Longrightarrow \
\alpha ( \rho ) \ = \ 
\begin{picture}(9,4)\thicklines
\put(0,-2){\line(0,1){3}} \put(0,-2){\line(1,0){10}}
\put(1,-2){\line(0,1){3}} 
\put(9,-2){\line(0,1){3}} 
\put(10,-2){\line(0,1){3}} 
\put(2,-1){\line(0,1){2}} \put(2,-1){\line(1,0){6}} 
\put(3,-1){\line(0,1){2}}
\put(6,-1){\line(0,1){2}}
\put(7,-1){\line(0,1){2}}
\put(8,-1){\line(0,1){2}} 
\put(4,0){\line(0,1){1}} \put(4,0){\line(1,0){1}} 
\put(5,0){\line(0,1){1}} 
\put(-0.3,1.5){0} \put(0.7,1.5){1} \put(1.7,1.5){2} 
\put(2.7,1.5){3}  \put(3.7,1.5){4} \put(4.7,1.5){5} 
\put(5.7,1.5){6}  \put(6.7,1.5){7} \put(7.7,1.5){8}
\put(8.7,1.5){9}  \put(9.7,1.5){10}
\end{picture}
\]

\vspace{20pt}

\noindent
Clearly, the definition of $\alpha ( \rho )$ could also be stated 
without referring to pictures. That is, the rule for assigning 
the singletons of $\rho$ to doubletons (in order to create 
$\alpha ( \rho )$) can be expressed in plain algebraic terms.
Indeed, for every 1-element block $\{ i \}$ of $\rho$, exactly 
one the following two possibilities (1) and (2) applies:

\noindent
(1) Either there is no 2-element block $B= \{ p,q \}$ of $\rho$ 
such that $p<i<q$. In this case $i$ is assigned to the doubleton 
$\{ 0, n+1 \}$ that was added to $\rho$.

\noindent
(2) Or there exist 2-element blocks $B= \{ p,q \}$ of $\rho$ 
such that $p<i<q$. Due to the fact that $\rho$ is non-crossing, 
among these blocks there has to exist one, $B_o = \{ p_o , q_o \}$,
which is nested inside all the others (we have $p<p_o$ and 
$q>q_o$ for every block $B = \{ p,q \}$, $B \neq B_o$, such that 
$p<i<q$). In this case the singleton $i$ is assigned to the 
doubleton $B_o$.

The construction of $\alpha ( \rho )$ described above defines 
a map
\begin{equation}  \label{eqn:6.3} 
\alpha : NC_{\leq 2} (n) \to 
NC( \, \{ 0,1, \ldots , n+1 \} \, ),
\end{equation}
the ``assign-singletons-to-doubletons'' map. It is easily checked
that the image of $\alpha$ is
\begin{equation}  \label{eqn:6.4} 
\left\{ \pi \in NC( \, \{ 0,1, \ldots , n+1 \} \, ) \
\begin{array}{cl}
\vline & \mbox{$0 \ecpi n+1$ and $\pi$ has} \\
\vline & \mbox{no 1-element blocks}
\end{array}  \right\} ,
\end{equation}
where the notation ``$0 \ecpi n+1$'' in (\ref{eqn:6.4}) is a 
shorthand for ``$0$ and $n+1$ belong to the same block of $\pi$''.
It is also immediate that the map $\alpha$ from (\ref{eqn:6.3})
is one-to-one. The map
\begin{equation}  \label{eqn:6.5} 
\beta :
\left\{ \pi \in NC( \, \{ 0,1, \ldots , n+1 \} \, ) \
\begin{array}{cl}
\vline & \mbox{$0 \ecpi n+1$ and $\pi$ has} \\
\vline & \mbox{no 1-element blocks}
\end{array}  \right\} \to NC_{\leq 2}(n)
\end{equation}
which is inverse to $\alpha$ is described as follows. Let 
$\pi = \{ A_1, \ldots , A_p \}$ be a partition from the set
(\ref{eqn:6.4}), and say that $A_1$ is the block of $\pi$
that contains $0$ and $n+1$. Then
\begin{equation}  \label{eqn:6.6} 
\beta ( \pi ) = \{ B_2, \ldots , B_p \} \cup
\bigl\{ \ \{ i \} \, \mid i \in \{ 1, \ldots , n \} \setminus
(B_2 \cup \cdots \cup B_p) \, \bigr\} ,
\end{equation}   
where for every $2 \leq i \leq p$ we denoted 
$B_i := \{ \min (A_i)$, 
$\max (A_i) \} \subseteq \{ 1, \ldots , n \}$. }
\end{remark}

$\ $

\noindent
{\bf Proof of Theorem \ref{thm:6.2}.} 
Fix $t > 0$ for which we will prove that (\ref{eqn:6.2}) holds. 
We will prove this equality by showing that the distributions on 
its two sides have the same $\eta$-series:
\begin{equation}  \label{eqn:6.7}
\eta_{\Phi ( \nu \boxplus \gamma_t )} =
\eta_{ \bB_t ( \Phi ( \nu ))} .
\end{equation}

We first observe that on both sides of (\ref{eqn:6.7}) we have
series in $\ncserk$ that are of the form
\begin{equation}  \label{eqn:6.8}
\Bigl( \sum_{i=1}^k z_i^2 \Bigr) + 
\mbox{(terms of order $\geq 3$)} .
\end{equation}
Indeed, from the definition of $\Phi$ in Equation (\ref{eqn:6.1})
it is clear that $\eta_{\Phi ( \sigma )}$ is of the form 
(\ref{eqn:6.8}) for every $\sigma \in \Dalg (k)$, and this applies
in particular to the left-hand side of (\ref{eqn:6.7}). On the 
right-hand side of (\ref{eqn:6.7}) we first invoke Remark 
\ref{rem:4.4} and write
\[
\eta_{ \bB_t ( \Phi ( \nu ))} =
\frac{1}{t} \Reta \Bigl( t \eta_{\Phi ( \nu )} \Bigr) ;
\]
then we use the fact that $\eta_{\Phi ( \nu )}$ is of the form 
(\ref{eqn:6.8}), combined with the observation (see Remark
\ref{rem:3.45}) that applying $\Reta$ does not change the linear
and quadratic terms of a series in $\ncserk$.

In order to prove (\ref{eqn:6.7}), we should thus fix a monomial 
of length $\geq 3$ in $z_1, \ldots , z_k$, and prove that the 
coefficients for this monomial in 
$\eta_{\Phi ( \nu \boxplus \gamma_t )}$ and in
$\eta_{ \bB_t ( \Phi ( \nu ))}$ are equal to each other. It will
be convenient to denote our fixed monomial in $z_1, \ldots , z_k$
as $z_{i_0} z_{i_1} \cdots z_{i_n} z_{i_{n+1}}$ for some $n \geq 1$
and $i_0, i_1, \ldots , i_{n+1} \leq k$. Our job for the 
remaining of the proof is to verify that 
\begin{equation}  \label{eqn:6.9}
\cf_{(i_0,i_1, \ldots , i_{n+1})} \Bigl( \,
\eta_{\Phi ( \nu \boxplus \gamma_t )} \, \Bigr) =
\cf_{(i_0,i_1, \ldots , i_{n+1})} \Bigl( \,
\eta_{ \bB_t ( \Phi ( \nu ))} \, \Bigr) ,
\end{equation}
for this fixed $n$ and $i_0, i_1, \ldots , i_{n+1}$.

On the left-hand side of (\ref{eqn:6.9}) we have
\begin{align*}
\cf_{(i_0,i_1, \ldots , i_{n+1})} \Bigl( \,
\eta_{\Phi ( \nu \boxplus \gamma_t )} \, \Bigr) 
& = \delta_{i_0, i_{n+1}} \cdot
    \cf_{(i_1, \ldots , i_n)} (M_{\nu \boxplus \gamma_t})
    \ \mbox{ (by Equation (\ref{eqn:6.1})) }               \\
& = \delta_{i_0, i_{n+1}} \cdot
    ( \nu \boxplus \gamma_t ) (X_{i_1} \cdots X_{i_n}).
\end{align*}
The latter moment is exactly of the kind studied in Section 5 of 
the paper, and can be expressed (by Proposition \ref{prop:5.5})
as a summation indexed by $NC_{\leq 2}(n)$. Thus for the 
left-hand side of (\ref{eqn:6.9}) we come to
\begin{equation}  \label{eqn:6.10}
\cf_{(i_0,i_1, \ldots , i_{n+1})} \Bigl( \,
\eta_{\Phi ( \nu \boxplus \gamma_t )} \, \Bigr) =
\delta_{i_0, i_{n+1}} \cdot \sum_{\rho \in NC_{\leq 2}(n)} \
\term_{\rho} ' ,
\end{equation}
where for every $\rho \in NC_{\leq 2}(n)$ the contribution 
$\term_{\rho} '$ of $\rho$ is as on the right-hand side of 
Equation (\ref{eqn:5.051}) in Proposition \ref{prop:5.5}.

On the right-hand side of Equation (\ref{eqn:6.9}) we go as
follows:
\[
\cf_{(i_0,i_1, \ldots , i_{n+1})} \Bigl( \,
\eta_{ \bB_t ( \Phi ( \nu ))} \, \Bigr) 
=  \cf_{(i_0,i_1, \ldots , i_{n+1})} \Bigl( \,
    \frac{1}{t} \Reta ( t \eta_{\Phi ( \nu )} ) \Bigr)
    \ \mbox{ (by Remark \ref{rem:4.4})}               
\]
\[
= \frac{1}{t} \sum_{\begin{array}{c}
{\scriptstyle \pi \in NC( \, \{ 0,1, \ldots , n+1 \} \, ) } \\
{\scriptstyle such \ that \ 0 \ecpi n+1}
\end{array}  } \ \cf_{(i_0,i_1, \ldots , i_{n+1}); \pi} 
    ( t \eta_{\Phi ( \nu )} )
    \ \mbox{ (by Proposition \ref{prop:3.4})}               \\
\]
\begin{equation}  \label{eqn:6.11}
= \sum_{\begin{array}{c}
{\scriptstyle \pi \in NC( \, \{ 0,1, \ldots , n+1 \} \, ) } \\
{\scriptstyle such \ that \ 0 \ecpi n+1}
\end{array}  } \ t^{| \pi |-1}
\cf_{(i_0,i_1, \ldots , i_{n+1}); \pi} ( \eta_{\Phi ( \nu )} ).
\end{equation}
Observe that the summation in (\ref{eqn:6.11}) may in fact be 
restricted to those partitions in 
$\pi \in NC( \, \{ 0,1, \ldots , n+1 \} \, )$ which (in 
addition to the condition that $0 \ecpi n+1$) are required 
to have no singleton blocks; this is because 
$\eta_{\Phi ( \nu ) }$ has no linear terms (see the discussion 
around (\ref{eqn:6.8}) above), thus
$\cf_{(i_0, i_1, \ldots , i_{n+1}); \pi} 
\bigl( \eta_{\Phi ( \nu ) } \bigr)= 0$ whenever $\pi$ has 
singleton blocks. So for the right-hand side of (\ref{eqn:6.9})
we arrive to the formula 
\begin{equation}  \label{eqn:6.12}
\cf_{(i_0,i_1, \ldots , i_{n+1})} \Bigl( \,
\eta_{ \bB_t ( \Phi ( \nu ))} \, \Bigr) 
= \sum_{\pi} \term_{\pi} '' ,
\end{equation}
where $\pi$ runs precisely in the set described in (\ref{eqn:6.4})
of Remark \ref{rem:6.3}, and where for such $\pi$ we put
\[
\term_{\pi} '' \ := \ t^{| \pi |-1} \cdot
\cf_{(i_0, i_1, \ldots , i_{n+1}); \pi}
\Bigl( \, \eta_{\Phi ( \nu )} \, \Bigr)      
\]
\begin{equation}  \label{eqn:6.115}
= t^{| \pi |-1} \cdot \prod_{ \begin{array}{c}
{\scriptstyle A \ block \ of \ \pi}               \\
{\scriptstyle A = \{ m_1 < m_2 < \cdots < m_p \} }
\end{array} } \ \delta_{i_{m_1}, i_{m_p}} \,
\nu \bigl( X_{i_{m_2}} \cdots X_{i_{m_{p-1}}} \bigr) .
\end{equation}
When writing (\ref{eqn:6.115}) we also took into account how 
$\eta_{\Phi ( \nu )}$ is defined by Equation (\ref{eqn:6.1}).

Let us next observe that if the indices 
$i_0, i_1, \ldots , i_{n+1}$ fixed since (\ref{eqn:6.9}) do not 
satisfy the condition $i_0 = i_{n+1}$, then the right-hand sides
of both (\ref{eqn:6.10}) and (\ref{eqn:6.12}) vanish. This is 
clear for (\ref{eqn:6.10}), while for (\ref{eqn:6.12}) we argue 
as follows: if $i_0 \neq i_{n+1}$ then the product in 
(\ref{eqn:6.115}) is guaranteed to vanish 
(since one of the blocks of $\pi$ contains 
$0$ and $n+1$), hence every term $\term_{\pi} ''$ on the 
right-hand side of (\ref{eqn:6.12}) is equal to $0$. 

So let us then assume that $i_0 = i_{n+1}$. The equality 
(\ref{eqn:6.9}) that we have to prove is reduced (by virtue of 
(\ref{eqn:6.10}) and (\ref{eqn:6.12})) to
\begin{equation}  \label{eqn:6.15}
\sum_{\rho \in NC_{\leq 2}(n)} \ \term_{\rho} ' \ = \ 
\sum_{ \begin{array}{c}
{\scriptstyle \pi \ in \ the}     \\
{\scriptstyle set \ from \ (\ref{eqn:6.4})} 
\end{array}  } \ \term_{\pi} '' ,
\end{equation}
where the quantities $\term_{\rho} '$ and $\term_{\pi} ''$ are 
described in Equations (\ref{eqn:5.051}) and (\ref{eqn:6.11}), 
respectively.
But the equality (\ref{eqn:6.15}) is immediately verified by 
using the ``assign-singletons-to-doubletons'' construction from
Remark \ref{rem:6.3}. Indeed, in Remark \ref{rem:6.3} we pointed
out a natural bijection $\beta$ from the set in (\ref{eqn:6.4})
onto $NC_{\leq 2}(n)$, and by using the explicit description 
provided there for $\beta$ it is immediately seen that 
$\term_{\pi} '' = \term_{\beta ( \pi )} '$,
for every $\pi$ in the set from (\ref{eqn:6.4}).
Thus $\beta$ provides a term-by-term identification of the sums 
on the two sides of (\ref{eqn:6.15}), and this completes the proof.
\hfill$\square$

$\ $

$\ $

\setcounter{section}{7}
\begin{center}
{\large\bf 7. Restricting to the framework of $\cD_c (k)$}
\end{center}

In this section we show that the results from the Sections 4--6 of
the paper continue to hold when we work in $C^*$-framework.

\setcounter{equation}{0}
\setcounter{theorem}{0}

\vspace{10pt}

\begin{definition}  \label{def:7.1}
{\rm 
We denote
\begin{equation}   \label{eqn:7.011}
\cD_c (k) = \left\{  \mu \in \Dalg (k)  \begin{array}{cc}
\vline  &  \exists \mbox{ $C^*$-probability space 
           $( \cA , \varphi )$}                          \\
\vline  &  \mbox{and selfadjoint elements 
          $x_1, \ldots , x_k \in \cA$}                   \\
\vline  &  \mbox{such that $\mu_{x_1, \ldots , x_k} = \mu$}
\end{array}  \right\} 
\end{equation}
(where the joint distribution $\mu_{x_1, \ldots , x_k}$ is 
defined as in Equation (\ref{eqn:3.102}) from Definition 
\ref{def:2.1}). The fact that $( \cA , \varphi )$ is a 
$C^*$-probability space means here that $\cA$ is a unital 
$C^*$-algebra and that $\varphi : \cA \to \bC$ is a positive 
linear functional such that $\varphi ( 1_{\cA} ) = 1$. 

The notation ``$\cD_c (k)$'' is chosen to remind of ``distributions 
with compact support'' -- indeed, in the case when $k=1$
we have a natural identification between $\cD_c (1)$ and the set of 
probability distributions with compact support on $\bR$.   }
\end{definition}

\vspace{10pt}

\begin{remark}  \label{rem:7.2}
{\rm 
In the preceding sections, the operations $\boxplus$ and $\uplus$ 
and the convolution powers with respect to them were considered
in the larger framework of the space $\Dalg (k)$. But by considering
sums of freely independent and respectively Boolean independent 
$k$-tuples of selfadjoint elements in a $C^*$-probability space, 
one sees that if $\mu , \nu \in \cD_c (k)$ then 
$\mu \boxplus \nu$ and $\mu \uplus \nu$ belong to $\cD_c (k)$
as well. Hence $\boxplus$ and $\uplus$ make sense as binary 
operations on $\cD_c (k)$. Moreover, concerning convolution 
powers we have that 
\begin{equation}  \label{eqn:7.021}
\mu \in \cD_c (k) \ \Rightarrow \ \left\{
\begin{array}{cll}
\mbox{(a)} & \mu^{\boxplus t} \in \cD_c (k) &
                             \forall \, t \geq 1, \mbox{ and }  \\
\mbox{(b)} & \mu^{\uplus t} \in \cD_c (k) & \forall \, t >0.
\end{array}  \right.  
\end{equation}
The fact stated in (\ref{eqn:7.021}(a)) was proved in
\cite{NS96}, by using compressions with free projections. 
The proof of (\ref{eqn:7.021}(b)) is done by constructing an 
operator model for $\mu^{\uplus t}$ -- see Remark 4.7 and 
Proposition 4.8 of \cite{BN06}. 

The following result is then an immediate consequence of 
(\ref{eqn:7.021}) and of what was proved in algebraic framework
in Theorem \ref{thm:4.11}.  }
\end{remark}

\begin{corollary}  \label{cor:7.25}

$1^o$ For every $t \geq 0$ it makes sense to define 
$\bB_t : \cD_c (k) \to \cD_c (k)$ by the formula
\begin{equation}  \label{eqn:7.031}
\bB_t ( \mu ) = \Bigl( \, \mu^{\boxplus (1+t)} \, 
\Bigr)^{\uplus (1/(1+t))}, \ \ \forall \, \mu \in \cD_c (k).
\end{equation}

$2^o$ The transformations of $\cD_c (k)$ defined by 
(\ref{eqn:7.031}) form a semigroup: 
$\bB_s \circ \bB_t = \bB_{s+t}$, $\forall \, s,t \geq 0$.

$3^o$ For $t=1$ we have $\bB_1 ( \mu ) = \bB ( \mu )$, 
$\forall \, \mu \in \cD_c (k)$, where $\bB$
is the multi-variable Boolean Bercovici-Pata bijection 
from Theorems 1 and 1' of the paper \cite{BN06}.
\hfill$\square$
\end{corollary}

$\ $

In the remaining part of this section we will show that the
above Theorem \ref{thm:6.2} also carries through to 
the $C^*$-framework. The main point that needs to be addressed
is that the map $\Phi : \Dalg (k) \to \Dalg (k)$ introduced in 
Definition \ref{def:6.1} sends $\cD_c (k)$ into itself. We will 
prove this via an ``operator model'' for $\Phi$, described in 
the next remark and theorem.  

$\ $

\begin{remark}  \label{rem:7.3}
(The operator model for $\Phi$.)
{\rm The {\em input} for this operator model is a system 
\[
( \cH ; a_1, \ldots , a_k ; \xi_o )
\]
where $\cH$ is a Hilbert space, $a_1, \ldots , a_k \in B( \cH )$
are selfadjoint operators, and $\xi_o \in \cH$ is a unit vector.
Starting from this data, we proceed as follows:

\vspace{10pt}

(i) We consider the Hilbert space
$\cK := \bC \oplus \left( \bigoplus_{j=1}^k \cH \right)$, and 
the unit vector $\Omega_0 :=1 \oplus
\underbrace{0\oplus0\oplus\cdots\oplus0}_{k\ {\rm times}}\in \cK$.
For $1 \leq j \leq k$ we let 
$v_j \colon \cH \to \cK$ be the embedding defined by
\[
v_j(\xi)= 0 \oplus
\underbrace{0 \oplus \cdots \oplus0}_{j-1\ {\rm times}}
\oplus \xi \oplus
\underbrace{0\oplus\cdots\oplus0}_{k-j\ {\rm times}} \in \cK , 
\ \mbox{ } \xi \in \cH .
\]
The direct sum defining $\cK$ can thus also be writtten as 
$\cK = \bC \Omega_0 \oplus v_1 ( \cH ) \oplus \cdots \oplus
v_k ( \cH )$.

\vspace{10pt}

(ii) For $1 \leq j \leq k$ we denote 
$v_j \xi_0=: \Omega_j \in \cK$, and we consider the 
rank-one partial isometry $w_j \in B( \cK )$ which carries 
$\Omega_0$ to $\Omega_j$. The operator $w_j$ and its adjoint 
are thus described by the formulas:
\[
w_j \eta =  \langle \eta \ , \ \Omega_0  \rangle \, \Omega_j , 
\ \mbox{ and } \
w_j^* \eta =  \langle \eta \ , \ \Omega_j  \rangle \, \Omega_0 , 
\ \ \forall \, \eta \in \cK .
\]

\vspace{10pt}

(iii) For $1 \leq j \leq k$ we consider the selfadjoint operators
$x_j , y_j \in B( \cK )$ defined by
\[
x_j := 0 \oplus
\underbrace{a_j \oplus \cdots \oplus a_j}_{k\ {\rm times}}
\ \mbox{ and } \ 
y_j := w_j+x_j+w_j^* .
\]

\vspace{10pt}

\noindent
The system $( \cK ; y_1, \ldots , y_k; \Omega_0 )$ will be called 
the {\em output} of the operator model for $\Phi$. 
The terms ``input'' and ``output'' used in this construction 
are justified by the following theorem.  }
\end{remark}

\vspace{10pt}

\begin{theorem}  \label{thm:7.4}
Let $( \cH ; a_1, \ldots , a_k; \xi_o )$ and 
$( \cK ; y_1, \ldots , y_k; \Omega_0 )$ be as in 
Remark \ref{rem:7.3}. Let $\nu$ be the joint distribution
of $a_1, \ldots , a_k$ with respect to the vector-state 
$\langle \ \cdot \ \xi_o  \ , \ \xi_o \rangle$ on $B( \cH )$,
and let $\mu$ be the joint distribution of $y_1, \ldots , y_k$ 
with respect to the vector-state 
$\langle \ \cdot \ \Omega_0  \ , \ \Omega_0 \rangle$ on $B( \cK )$.
Then $\Phi ( \nu ) = \mu$.
\end{theorem}

\vspace{10pt}

\begin{remark}  \label{rem:7.5}
{\rm In preparation of the proof of Theorem \ref{thm:7.4} we
review here the ``moments vs. Boolean cumulants'' formula, which 
expresses the moments of a distribution $\mu \in \Dalg (k)$ 
in terms of its Boolean cumulants -- that is, in terms of the 
coefficients of the $\eta$-series $\eta_{\mu}$. This is very 
similar to the moment-cumulant formula reviewed in Remark 
\ref{rem:5.7} in connection to free cumulants,
with the difference that we now only consider summations over 
the subposet of $NC(n)$ consisting of interval partitions.

A partition $\pi$ of $\{ 1, \ldots , n \}$ is said to be an 
{\em interval partition} when every block of $\pi$ is of the form 
$[ a,b ] \cap \bZ$ for some $a \leq b$ in $\{ 1, \ldots , n \}$.
The set of all interval partitions of $\{ 1, \ldots , n \}$ will
be denoted by $\Int (n)$. It is clear that 
$\Int (n) \subseteq NC(n)$. The ``momemnts vs. Boolean cumulants''
formula says that for a distribution $\mu \in \Dalg (k)$ we have
\begin{equation}  \label{eqn:7.051}
\cf_{(i_1, \ldots , i_n)} (M_{\mu}) = \sum_{\pi \in \Int (n)} 
\cf_{(i_1, \ldots , i_n); \pi} ( \eta_{\mu} ), 
\ \ \forall \, n \geq 1, \ \forall \,
1 \leq i_1, \ldots , i_n \leq k.
\end{equation}
Equation (\ref{eqn:7.051}) is easily seen to be equivalent to the 
formula ``$\eta_{\mu} = M_{\mu} \bigl( 1+ M_{\mu} \bigr)^{-1}$''
used in the above Remark \ref{rem:2.4} as definition for the 
$\eta$-series of $\mu$ (for a proof of this equivalence, see 
for instance Proposition 3.5 in \cite{BN06}).  }
\end{remark}

\vspace{10pt}

\begin{remark}  \label{rem:7.6}
{\rm We now return to the notations from Remark \ref{rem:7.3},
and record how the operators $w_i, x_i, w_i^*$ 
$( 1 \leq i \leq k)$ behave with respect to the direct sum 
decomposition $\cK = \bC \Omega_0 \oplus v_1 ( \cH ) 
\oplus \cdots \oplus v_k ( \cH )$: we have that 
\begin{equation}  \label{eqn:7.061}
\left\{  \begin{array}{l}
\mbox{$w_i$ sends $\bC \Omega_0$ to $v_i ( \cH )$ and sends
      $v_1 ( \cH ), \ldots , v_k ( \cH )$ to $0$;}           \\
                                                             \\
\mbox{$x_i$ sends $\bC \Omega_0$ to $0$ and sends every 
      $v_i ( \cH )$ into itself, $1 \leq i \leq k$;}         \\
                                                             \\
\mbox{$w_i^*$ sends $v_i ( \cH )$ into $\bC \Omega_0$ and 
      sends $\bC \Omega_0$ and every $v_j ( \cH )$ with 
      $j \neq i$ to $0$.}           
\end{array}  \right.
\end{equation}
The verification of (\ref{eqn:7.061}) is immediate from the 
explicit formulas describing $w_i, x_i, w_i^*$ in Remark
\ref{rem:7.3}.    }
\end{remark}

\vspace{10pt}

\begin{lemma}  \label{lemma:7.7}
Consider the notations from Remark \ref{rem:7.3}, and let $\nu$
denote the joint distribution of $a_1, \ldots , a_k$ with respect 
to the vector-state $\langle \ \cdot \ \xi_o \ , \ \xi_o \rangle$ 
on $B( \cH )$. Let $j_1, \ldots , j_m$ and $i',i''$ be some indices 
in $\{ 1, \ldots , k \}$. Then we have
\begin{equation}  \label{eqn:7.071}
w_{i'}^* x_{j_1} \cdots x_{j_m} w_{i''} \Omega_0
= \lambda \Omega_0 ,
\end{equation}
where $\lambda = \cf_{(i',j_1, \ldots ,j_m,i'')}
\Bigl( \, \eta_{\Phi ( \nu )} \, \Bigr)$.
\end{lemma}

$\ $

\noindent
{\bf Proof.} If $i' \neq i''$ then both sides of Equation 
(\ref{eqn:7.071}) are equal to 0: the right-hand side vanishes
because of how $\eta_{\Phi ( \nu )}$ is defined (see Definition
\ref{def:6.1}), while the vanishing on the left-hand side follows
immediately from the operating rules described in 
(\ref{eqn:7.061}).  So we will assume that $i' = i'' =: i$, when 
the relation that has to be proved becomes;
\[
w_{i}^* x_{j_1} \cdots x_{j_m} w_{i} \Omega_0
= \nu ( X_{j_1} \cdots X_{j_m} )  \Omega_0 .
\]
We have $w_i ( \Omega_0 ) = \Omega_i$, and directly from the 
definition of $x_1, \ldots , x_k$ we observe that
$x_{j_1} \cdots x_{j_m} \Omega_i = v_i 
(a_{j_1} \cdots a_{j_m}  \xi_o )$.
But then:
\begin{align*}
w_{i}^* x_{j_1} \cdots x_{j_m} w_{i} \Omega_0
& = w_i^* v_i (a_{j_1} \cdots a_{j_m}  \xi_o )           \\
& = \langle v_i (a_{j_1} \cdots a_{j_m}) \xi_o \ , 
    \ \Omega_i \rangle \Omega_0                          \\
& = \langle a_{j_1} \cdots a_{j_m} \xi_o \ , 
    \ \xi_o \rangle \Omega_0 \ \ \mbox{ (since 
    $v_i^* \Omega_i = \xi_o$)}                           \\
& = \nu ( X_{j_1} \cdots X_{j_m} )  \Omega_0 ,
\end{align*}
as required.
\hfill$\square$

$\ $

\begin{lemma}  \label{lemma:7.8}
Consider the notations from Remark \ref{rem:7.3}, and let $\nu$
denote the joint distribution of $a_1, \ldots , a_k$ with respect 
to the vector-state $\langle \ \cdot \ \xi_o \ , \ \xi_o \rangle$ 
on $B( \cH )$. Let $i_1, \ldots , i_n$ be some indices in 
$\{ 1, \ldots , k \}$. Let $\pi$ be a partition in 
$\Int (n)$ which has no 1-element blocks, and which is written
explicitly as
$\pi = \{ \ \{ a_1, \ldots , b_1 \}, \ldots ,
\{ a_p, \ldots , b_p \} \ \}$, with 

\noindent
$1 = a_1 < b_1 < \cdots < a_p < b_p = n$ (and where 
$a_2 = b_1 + 1, \ldots , a_p = b_{p-1} +1$).
Consider the operators $u_1, \ldots , u_n \in B( \cK )$ defined
as follows:
\begin{equation}  \label{eqn:7.081}
\left\{  \begin{array}{l}
u_{a_1} = w_{i_{a_1}}^* , \ldots , u_{a_p} = w_{i_{a_p}}^*, \\
                                                    \\
u_{b_1} = w_{i_{b_1}}   , \ldots , u_{b_p} = w_{i_{b_p}},   \\
                                                    \\
u_c = x_{i_c} \mbox{ for every }
  c \in \{ 1, \ldots , n \} \setminus \{ a_1, b_1, \ldots ,
                                       a_p, b_p \} .
\end{array}  \right.
\end{equation}
Then we have
\begin{equation}  \label{eqn:7.082}
\langle u_1 \cdots u_n \Omega_0 \ , \ \Omega_0 \rangle 
= \cf_{( i_1,  \ldots , i_n ); \pi} 
\Bigl( \, \eta_{\Phi ( \nu )} \, \Bigr) .
\end{equation}
\end{lemma}

$\ $

\noindent
{\bf Proof.} By picking out the last $b_p - a_p +1$ factors in the 
product $u_1 \cdots u_n$ applied to the vector $\Omega_0$ we get:
\begin{equation}  \label{eqn:7.083}
u_{a_p} u_{a_p + 1} \cdots u_{b_p -1} u_{b_p} \Omega_0
= w_{i_{a_p}}^* \cdot \prod_{a_p < c < b_p} x_{i_c} \cdot
w_{i_{b_p}} \Omega_0
= \cf_{( i_{a_p}, i_{a_p +1}, \ldots , i_{b_p} )} \
\Bigl( \, \eta_{\Phi ( \nu )} \, \Bigr) \Omega_0,
\end{equation}  
where at the second equality sign we invoked Lemma 
\ref{lemma:7.7}. Thus
\begin{align*}
u_1 \cdots u_n \Omega_0 
& = u_1 \cdots u_{b_{p-1}} \bigl( 
u_{a_p} \cdots u_{b_p} \Omega_0 \bigr)               \\
& = \cf_{( i_{a_p}, i_{a_p +1}, \ldots , i_{b_p} )} \
\Bigl( \, \eta_{\Phi ( \nu )} \, \Bigr) \cdot
( u_1 \cdots u_{b_{p-1}} \Omega_0).
\end{align*}
Now the same trick as in (\ref{eqn:7.083}) can be 
applied to the right-most piece  
$u_{a_{p-1}} \cdots u_{b_{p-1}} \Omega_0$ of 
$u_1 \cdots u_{b_{p-1}} \Omega_0$. By iterating this trick 
we arrive to required the conclusion that
\begin{align*}
\langle u_1 \cdots u_n \Omega_0 \ , \ \Omega_0 \rangle
& = \cf_{( i_{a_1}, \ldots , i_{b_1} )} 
\Bigl( \, \eta_{\Phi ( \nu )} \, \Bigr) \cdots
\cf_{( i_{a_p}, \ldots , i_{b_p} )} 
\Bigl( \, \eta_{\Phi ( \nu )} \, \Bigr)           \\ 
& = \cf_{( i_1,  \ldots , i_n ); \pi} 
\Bigl( \, \eta_{\Phi ( \nu )} \, \Bigr) .
\hspace{2cm} \square
\end{align*}

$\ $

\noindent
{\bf Proof of Theorem \ref{thm:7.4}.}
We fix for the whole proof a positive integer $n$ and some
indices $1 \leq i_1, \ldots , i_n \leq k$, for which we verify 
that
\begin{equation} \label{eqn:7.4}
\mu ( X_{i_1} \cdots X_{i_n} ) 
= \bigl( \Phi ( \nu ) \bigr) ( X_{i_1} \cdots X_{i_n} ).
\end{equation}
By the definition of $\mu$, the left-hand side of (\ref{eqn:7.4})
is
\begin{align*}
\mu ( X_{i_1} \cdots X_{i_n} ) 
& = \langle y_{i_1} \cdots y_{i_n} \Omega_0 \ , 
    \ \Omega_0 \rangle                                \\
& = \langle (w_{i_1} + x_{i_1} + w_{i_1}^*) \cdots 
     (w_{i_1} + x_{i_1} + w_{i_1}^*) \Omega_0 \ , 
    \ \Omega_0 \rangle ,
\end{align*}
and the latter quantity expands as a sum of $3^n$ terms of the
form $\langle u_1 \cdots u_n \Omega_0 \ , \ \Omega_0 \rangle$,
with
\begin{equation}  \label{eqn:7.5}
u_1 \in \{ w_{i_1}, x_{i_1}, w_{i_1}^* \}, \ldots ,
             u_n \in \{ w_{i_n}, x_{i_n}, w_{i_n}^* \}.
\end{equation}
But from the rules (\ref{eqn:7.061}) for how 
the operators
$w_i, x_i, w_i^*$ act on the decomposition 
$\bC \Omega_0 \oplus v_1 ( \cH ) \oplus \cdots \oplus 
v_k ( \cH )$ of $\cK$ it follows that many of these $3^n$ terms
vanish. We leave it as an easy exercise to the reader to verify
that we have in fact 
$\langle u_1 \cdots u_n \Omega_0 \ , \ \Omega_0 \rangle$ = 0
whenever $u_1, \ldots , u_n$ from (\ref{eqn:7.5}) are not chosen
according to the recipe (\ref{eqn:7.081}) from Lemma 
\ref{lemma:7.8}. By taking Lemma \ref{lemma:7.8} into account,
we thus find that
\begin{equation}   \label{eqn:7.11}
\mu (X_{i_1} \cdots X_{i_n}) = \sum_{\begin{array}{c}
{\scriptstyle \pi \in \Int (n), \ with}  \\
{\scriptstyle no \ singletons}  
\end{array}  } \
\cf_{( i_1,  \ldots , i_n ); \pi} 
\Bigl( \, \eta_{\Phi ( \nu )} \, \Bigr) .
\end{equation}

It remains to note that on the right-hand side of 
(\ref{eqn:7.11}) it does not cost anything to add the terms
$\cf_{( i_1,  \ldots , i_n ); \pi} 
\Bigl( \, \eta_{\Phi ( \nu )} \, \Bigr)$ where $\pi \in \Int (n)$
has some singleton blocks; indeed, each of these added terms is 
in fact equal to $0$, because the linear terms of the series 
$\eta_{\Phi ( \nu )}$ vanish. So from (\ref{eqn:7.11}) we can write:
\begin{align*}
\mu (X_{i_1} \cdots X_{i_n}) 
& = \sum_{\pi \in \Int (n)} \cf_{( i_1,  \ldots , i_n ); \pi} 
    \Bigl( \, \eta_{\Phi ( \nu )} \, \Bigr)                    \\
& = \cf_{( i_1,  \ldots , i_n )} \Bigl( \, M_{\Phi ( \nu )} 
    \, \Bigr) \ \ \mbox{ (by Remark \ref{rem:7.5})}            \\
& = \Bigl( \, \Phi ( \nu ) \, \Bigr) (X_{i_1} \cdots X_{i_n}), 
\end{align*}
which is what we wanted to obtain.
\hfill$\square$

$\ $

\begin{corollary}  \label{cor:7.9}
The map $\Phi : \Dalg (k) \to \Dalg (k)$ introduced in 
Definition \ref{def:6.1} carries $\cD_c (k)$ into itself.
\end{corollary}

\vspace{6pt}

\noindent
{\bf Proof.} Let $\nu$ be in $\cD_c (k)$. By using the GNS
construction one can realize $\nu$ as the joint distribution 
of a $k$-tuple $a_1, \ldots , a_k$ of selfadjoint operators 
on a Hilbert space $\cH$, with respect to a vector-state 
$\langle \ \cdot \ \xi_o \ , \ \xi_o \rangle$ on $B( \cH )$. 
Then Theorem \ref{thm:7.4} gives $\Phi ( \nu )$ as the joint 
distribution of $y_1, \ldots , y_k \in B( \cK )$ with 
respect to $\langle \ \cdot \ \Omega_0 \ , \ \Omega_0 \rangle$,
where $( \cK ; y_1, \ldots , y_k ; \Omega_0 )$ are constructed
as in Remark \ref{rem:7.3}. This implies that 
$\Phi ( \nu ) \in \cD_c (k)$.
\hfill $\square$

$\ $

It thus follows that the statement of Theorem \ref{thm:6.2} 
also holds in $C^*$-framework:

\begin{corollary}  \label{cor:7.10}
Let $\nu$ be a distribution in $\cD_c (k)$. Then for 
every $t>0$ we have
\begin{equation}  \label{eqn:7.101}
\Phi ( \, \nu \boxplus \gamma_t \, ) = 
\bB_t ( \, \Phi ( \nu ) \, ) \in \cD_c (k),
\end{equation}
where $\gamma_t$ is the distribution of the scaled free 
semicircular system from Notation \ref{def:5.1}.
\hfill $\square$
\end{corollary}

$\ $

$\ $

\newpage

Serban T. Belinschi: University of Waterloo and IMAR.

Address: Department of Pure Mathematics, University of Waterloo,  

Waterloo, Ontario N2L 3G1, Canada.                      

Email: sbelinsc@math.uwaterloo.ca

$\ $

Alexandru Nica: University of Waterloo.

Address: Department of Pure Mathematics, University of Waterloo, 

Waterloo, Ontario N2L 3G1, Canada.                    

Email: anica@math.uwaterloo.ca

\end{document}